\documentclass[reprint, amsmath,amssymb,aps,pra]{revtex4-2}

\usepackage[pdfborder={0 0 0.5 [3 2]}, plainpages=false]{hyperref}%
\usepackage[left=1in,right=1in,top=1in,bottom=1in]{geometry}%
\usepackage{amsmath}
\usepackage{enumerate}
\usepackage{amssymb}                
\usepackage{amsmath}                
\usepackage{amsfonts}
\usepackage{amsthm}
\usepackage{bbm}
\usepackage[table,xcdraw]{xcolor}
\usepackage{booktabs}
\usepackage{svg}
\usepackage{mathtools}
\usepackage{cool}
\usepackage{url}
\usepackage{graphicx,epsfig}
\usepackage{makecell}
\usepackage{array}

\usepackage[capitalize,nameinlink]{cleveref}
\crefname{section}{section}{sections}
\crefname{subsection}{subsection}{subsections}
\Crefname{section}{Section}{Sections}
\Crefname{subsection}{Subsection}{Subsections}

\Crefname{figure}{Figure}{Figures}

\crefformat{equation}{\textup{#2(#1)#3}}
\crefrangeformat{equation}{\textup{#3(#1)#4--#5(#2)#6}}
\crefmultiformat{equation}{\textup{#2(#1)#3}}{ and \textup{#2(#1)#3}}
{, \textup{#2(#1)#3}}{, and \textup{#2(#1)#3}}
\crefrangemultiformat{equation}{\textup{#3(#1)#4--#5(#2)#6}}%
{ and \textup{#3(#1)#4--#5(#2)#6}}{, \textup{#3(#1)#4--#5(#2)#6}}{, and \textup{#3(#1)#4--#5(#2)#6}}

\Crefformat{equation}{#2Equation~\textup{(#1)}#3}
\Crefrangeformat{equation}{Equations~\textup{#3(#1)#4--#5(#2)#6}}
\Crefmultiformat{equation}{Equations~\textup{#2(#1)#3}}{ and \textup{#2(#1)#3}}
{, \textup{#2(#1)#3}}{, and \textup{#2(#1)#3}}
\Crefrangemultiformat{equation}{Equations~\textup{#3(#1)#4--#5(#2)#6}}%
{ and \textup{#3(#1)#4--#5(#2)#6}}{, \textup{#3(#1)#4--#5(#2)#6}}{, and \textup{#3(#1)#4--#5(#2)#6}}

\crefdefaultlabelformat{#2\textup{#1}#3}

\def\noi{\noindent}

\def\R{{\mathbb R}}

\DeclareMathOperator{\diag}{diag}

\graphicspath{ {images/} }

\begin{document}

\title{Standing Wave Solutions in Twisted Multicore Fibers}

\author{Ross Parker}
\address{Department of Mathematics, Southern Methodist University, 
Dallas, TX 75275}
\email{rhparker@smu.edu}

\author{Alejandro Aceves}
\address{Department of Mathematics, Southern Methodist University, 
Dallas, TX 75275}
\email{aaceves@smu.edu}

\begin{abstract}
In the present work, we consider the existence and spectral stability of standing wave solutions to a model for light propagation in a twisted multi-core fiber with no gain or loss of energy. Numerical parameter continuation experiments demonstrate the existence of standing wave solutions for sufficiently small values of the coupling parameter. Furthermore, standing waves exhibiting optical Aharonov-Bohm suppression, where there is a single waveguide which remains unexcited, exist when the twist parameter $\phi$ and the number of waveguides $N$ is related by $\phi = \pi/N$. Spectral computations and numerical evolution simulations suggest that standing wave solutions where the energy is concentrated in a single site are neutrally stable. Solutions with asymmetric coupling and multi-pulse solutions are also briefly explored. 
\end{abstract}

\maketitle

\section{Introduction}

There has been much recent theoretical and experimental interest in light dynamics in twisted multi-core optical fibers. Early work on twisted fibers can be found in \cite{Longhi2007,Longhi2007b}, in which the coupled mode equations describing light propagation in a twisted, circular arrangement of waveguides is derived. The introduction of a fiber twist in a circular array allows for control of diffraction and light transfer, in a similar manner to axis bending in linear waveguide arrays \cite{Longhi2005}. The fiber twist introduces additional phase terms to the model, which is known as the Peierls phase \cite{Longhi2007,Peierls1933}. In \cite{Ornigotti2007}, this system is considered as an optical analogue of topological Aharonov-Bohm suppression of tunneling \cite{Loss1992}, where the fiber twist plays the role of the magnetic flux in the quantum mechanical system. In the optical setting, what this suppression we believe reveals is similar to bending of rays in twisted photonic crystals \cite{Russell}, resulting in the creation of ``forbidden'' access points in the transverse profile as rays propagate in the longitudinal direction. Alternatively, the phase accumulation from the twist, together with that due to the amplitude-dependent phase differences, accounts for a phase mismatch that inhibits transfer of energy among waveguides. The unique feature present here is that the suppression is full, thus instead of a localized mode with nonzero amplitudes across the array, a topological state is achieved. This state is both nonlinear and robust. Fiber arrangements featuring parity-time ($\mathcal{PT}$) symmetry with balanced gain and loss terms are considered in \cite{Longhi2016,castro2016}. More complicated fiber bundle geometries have since been studied, which include Lieb lattices \cite{Marzuola2019bulk} and honeycomb lattices \cite{Ablowitz2014,Lumer2013}. Experimental applications of twisted multi-core fibers include the construction of sensors for shape, strain, and temperature \cite{Gannot2014,Westbrook2017}. 

\begin{figure}
\begin{center}
\includegraphics[width=3.5cm]{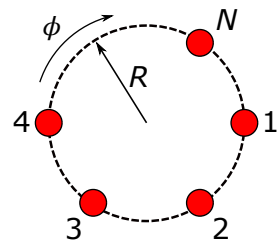}
\end{center}
\caption{Schematic of twisted, multi-core fiber consisting of $N$ waveguides arranged in a ring.}
\label{fig:ring}
\end{figure}

In this paper, we consider a multi-core fiber consisting of $N$ waveguides arranged in a ring (\cref{fig:ring}).
\noi The entire fiber is twisted in a uniform fashion along the propagation direction $z$ with a twist rate $\epsilon = 2\pi/\Lambda$, where $\Lambda$ is the spatial period. For the system with an optical Kerr nonlinearity, the dynamics are given by the coupled system of equations
\begin{equation}\label{eq:twist}
i \partial_z c_n = k \left(e^{-i\phi}c_{n+1} + e^{i\phi}c_{n-1}\right) + i \gamma_n c_n + d |c_n|^2 c_n
\end{equation}
for $n = 1, \dots, N$, where $c_0 = c_{N}$ and $c_{N+1} = c_1$ due to the circular geometry \cite{castro2016,Parto2017}. In the discrete approximation, the assumption is that the energy of electromagnetic field propagating along the optical array is concentrated in the guiding (silica) cores. As such, the complex amplitudes $c_n(z)$ represent the localized field amplitude in each waveguide. Since the tail of the transverse field profile at each waveguide extends beyond the core, the tail field concentrated at site $n$ overlaps with its neighbor cores at sites $n\pm1$. In this approximation, $k$ (in $\text{mm}^{-1}$ units) is the strength of the nearest-neighbor-waveguide coupling, $d$ is the effective (and normalized) Kerr-nonlinear index of refraction, and $\gamma_n$ is the optical gain (due to doping) or loss (due to imperfections or scattering) at waveguide $n$. Altogether, in the discrete approximation, all coefficients depend on the wavelength $\lambda$ and $\phi = 4 \pi^2 \epsilon n_s R^2/N \lambda$ is the Peierls phase introduced by the twist, where $n_s$ the refractive index of the substrate, $R$ is the radius of the circular ring, and $\lambda$ is the wavelength of the propagating field \cite{castro2016} (see also \cite{Longhi2007b,Garanovich2012} for a derivation of this equation). If $\gamma_n = 0$ for all $n$, i.e. there is no gain or loss at each node, the system is conservative. Furthermore, upon normalizing the fields by making $c_n$ non-dimensional using the mapping $c_n \mapsto \frac{1}{\sqrt{|d|}}c_n$,  equation \cref{eq:twist} becomes
\begin{equation}\label{eq:twist1}
i \partial_z c_n = k \left(e^{-i\phi}c_{n+1} + e^{i\phi}c_{n-1}\right)  \pm |c_n|^2 c_n,
\end{equation}
which is Hamiltonian with conserved energy
\begin{equation}\label{eq:H}
H = \sum_{n=1}^N k (c_{n+1}c_n^* e^{-i \phi} + c_n c_{n+1}^* e^{i \phi}) \pm \frac{1}{2}|c_n|^4.
\end{equation}
In this paper, we will only be concerned with the Hamiltonian system \cref{eq:twist1} with conserved quantity \cref{eq:H}; we will also only consider the defocusing (minus) nonlinearity. The case with symmetric gain-loss terms ($\mathcal{PT}$ symmetry) is considered in \cite{castro2016}. Asymptotic analysis of the system \cref{eq:twist1} for $N=6$ fibers where the peak intensity is contained in the first fiber ($n=1$) shows that the opposite fiber in the ring ($n=4$) has, to leading order, zero intensity when the twist parameter is given by $\phi = \pi/6$ \cite{castro2016}. This is confirmed by numerical time evolution simulations (see \cite[Figures 4 and 5]{castro2016}). This phenomenon is discussed in the context of Aharonov-Bohm (AB) suppression of optical tunneling in twisted multicore fibers in \cite{Parto2017,Parto2019}. In particular, this effect is demonstrated analytically for the case of $N = 4$ fibers and $\phi = \pi/4$ by solving the nonlinear system \cref{eq:twist} analytically \cite{Parto2019}. A natural question is whether AB suppression is present for larger N and if this state is robust (stable).

In what follows, we study the existence and stability of standing wave solutions (bound states) of equation \cref{eq:twist1}. This paper is organized as follows. In \cref{sec:standingwave}, we use numerical parameter continuation to construct standing wave solutions to \cref{eq:twist1} where the bulk of the energy is confined to a single fiber. In \cref{sec:ABsupp}, we demonstrate the existence, both analytically and numerically, of standing wave solutions which have a single dark node; this occurs when $\phi = \pi/N$, both for $N$ even and $N$ odd. We then investigate the stability of these solutions in \cref{sec:stability}. We conclude with a brief discussion of asymmetric variants and multi-modal solutions and suggest some directions for future research.

\section{Standing wave solutions}\label{sec:standingwave}

Standing wave solutions to \cref{eq:twist1} are bound states of the form
\begin{equation}\label{eq:ansatz1}
c_n = a_n e^{i (\omega z + \theta_n) },
\end{equation}
where $a_n \in \R$, $\theta_n \in (-\pi/2, \pi/2]$, and $\omega$ is the propagation constant. (Since we allow $a_n$ to be negative, we can restrict $\theta_n$ to that interval). We will refer to the $a_n$ as the amplitudes and the $\theta_n$ as the phases of each node. Standing waves are periodic in $z$ with period $2\pi/\omega$, and the intensity at each node $|c_n| = |a_n|$ is constant in $z$.
Making this substitution and simplifying, equation \cref{eq:twist1} becomes
\begin{equation}\label{eq:twisteq}
\begin{aligned}
k\Big( &a_{n+1} e^{i((\theta_{n+1}-\theta_n)-\phi)} \\
&+ a_{n-1} e^{-i((\theta_n - \theta_{n-1})-\phi)}\Big) + \omega a_n - a_n^3 = 0,
\end{aligned}
\end{equation}
where we have taken the defocusing (minus) nonlinearity. Equation \cref{eq:twisteq} can be written as the system of $2n$ equations
\begin{equation}\label{eq:twisteqreal}
\begin{aligned}
k\big( &a_{n+1} \cos(\theta_{n+1}-\theta_n-\phi) \\
&\quad+a_{n-1} \cos(\theta_n - \theta_{n-1}-\phi)\big) \\
&\quad+ \omega a_n -  a_n^3 = 0 \\
&a_{n+1} \sin(\theta_{n+1}-\theta_n-\phi) \\
&\quad-a_{n-1} \sin(\theta_n - \theta_{n-1}-\phi) = 0
\end{aligned}
\end{equation}
by separating real and imaginary parts. We note that the the exponential terms in \cref{eq:twisteq} depend only on the phase differences $\theta_{n+1}-\theta_n$ between adjacent sites. Due to the gauge invariance of \cref{eq:twist1}, if $c_n$ is solution, so is $e^{i \theta} c_n$, thus we may without loss of generality take $\theta_1 = 0$. If $\phi = 0$, i.e. the fibers are not twisted, we can take $\theta_n = 0$ for all $n$, and so \cref{eq:twisteq} reduces to the untwisted case with periodic boundary conditions. Similarly, if we take $\phi = 2 \pi/N$ and $\theta_n = (n-1)\phi$ for all $n$, the exponential terms do not contribute, and \cref{eq:twisteq} once again reduces to untwisted case. The interesting cases, therefore, occur when $0 < \phi < 2 \pi/N$. 

\begin{figure}
\begin{center}
\includegraphics[width=5cm]{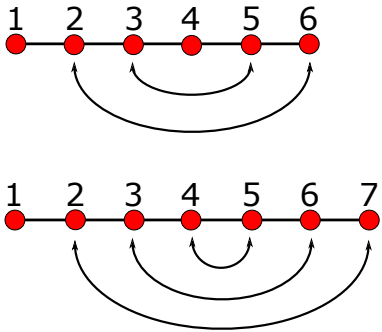}
\end{center}
\caption{Schematic of symmetry relationship between nodes for $N = 6$ and $N=7$. For nodes connected with arrows, the amplitudes $a_k$ are the same and the phases $\theta_k$ are opposite.}
\label{fig:symmetry1}
\end{figure}

In the anti-continuum (AC) limit ($k = 0$), the lattice sites are decoupled. Each $a_n$ can take on the values $\{0, \pm \sqrt{\omega} \}$, the phases $\theta_n$ are arbitrary, and $\phi$ does not contribute. The amplitudes $\sqrt{\omega}$ are real if  $\omega > 0$. We construct solutions to \cref{eq:twisteqreal} by parameter continuation from the AC limit with no twist using the standard continuation software package AUTO \cite{AUTO}. As an initial condition, we choose a single excited site at node 1, i.e. $a_1 = \sqrt{\omega}$ and $a_n = 0$ for all other $n$. (We can start with more than once excited state, but, in general, these solutions will not be stable.) In addition, we take $\theta_n = 0$ for all $n$, and $\phi = 0$.  We first continue in the coupling parameter $k$, and then, for fixed $k$, we continue in the twist parameter $\phi$. In doing this, we observe that the solutions have the following symmetry:
\begin{equation}\label{eq:symm}
\begin{aligned}
a_j &= a_{N-j+2} && \qquad j = 2, \dots, M-1 \\
\theta_j &= -\theta_{N-j+2} && \qquad j = 2, \dots, M-1,
\end{aligned}
\end{equation}
where $M = (N/2)+1$ for $N$ even and $M = (N+1)/2$ for $N$ odd. For $N$ even, node $M$ is the node directly across the ring from node 1, and $\theta_M = 0$. For all $N$, $\theta_1 = 0$. See \cref{fig:symmetry1} for an illustration of these symmetry relations for $N = 6$ and $N = 7$. 

\begin{figure}
\begin{center}
\begin{tabular}{c}
\includegraphics[width=8cm]{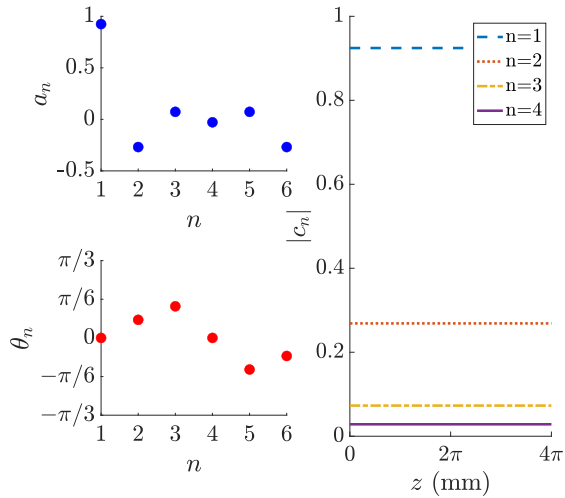}
\end{tabular}
\end{center}
\caption{Standing wave solution for $N = 6$, $\omega = 1 \text{ mm}^{-1}$, $k = 0.25 \text{ mm}^{-1}$, and $\phi = 0.25$. Left panel shows amplitudes $a_n$ and phases $\phi_n$ for solution at each node. Right panel is intensity of solution $|c_n|$ versus $z$ for nodes 1-4, which is constant in $z$. Evolution in $z$ computed using fourth order Runge-Kutta method.}
\label{fig:twist025}
\end{figure}

\cref{fig:twist025} shows an example of a standing wave solution of the form \cref{eq:ansatz1} produced by numerical parameter continuation for $N = 6$, $k = 0.25$, and $\phi = 0.25$. Since the paramater continuaton was initialized with a single excited site at node 1 in the AC limit, the peak intensity is still contained in node 1 when $k>0$, although the intensity has spread to the other nodes in the ring. The symmetry relations \cref{eq:symm} among the amplitudes $a_n$ and phases $\theta_n$ can be seen in the left panel. The node with minimum intensity is the node directly across the ring from node 1. The right panel shows the intensity $|c_n|$ at each node as a function of $z$. Since these are standing wave solutions, the intensity is constant in $z$. The evolution in $z$ is computed with a fourth-order Runga-Kutta method using equation \cref{eq:ansatz1} with $z=0$ and the amplitudes and phases from the left panel of \cref{fig:twist025} as the initial condition. This initial condition is used for all evolution plots for standing waves.

\begin{figure}
    \begin{center}
    \begin{tabular}{c}
    \includegraphics[width=8cm]{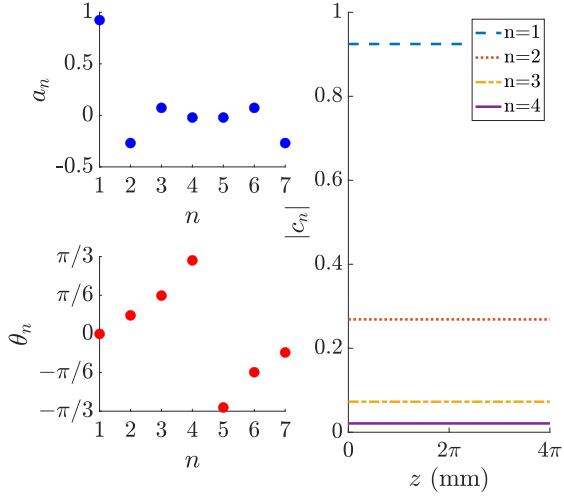}
    \end{tabular}
    \end{center}
    \caption{Standing wave solution for $N = 7$, $\omega = 1\text{ mm}^{-1}$, $k = 0.25\text{ mm}^{-1}$, and $\phi = 0.25$. Left panel shows amplitudes $a_n$ and phases $\phi_n$ for solution at each node. Right panel is intensity of solution $|c_n|$ versus $z$ for nodes 1-4, which is constant in $z$. Evolution in $z$ computed using fourth order Runge-Kutta method.}
    \label{fig:twist025N7}
    \end{figure}

Similarly, \cref{fig:twist025N7} shows a standing wave solution produced by numerical parameter continuation for $N = 7$, $k = 0.25$, and $\phi = 0.25$. As with the case of $N=6$, the peak intensity is contained in the node 1, and the symmetry relations \cref{eq:symm} among the amplitudes $a_n$ and phases $\theta_n$ can be seen in the left panel. In contrast with the $N=6$ case, there is a pair of nodes with minimum intensity and the same amplitude directly across the ring from node 1.

\begin{figure}
    \begin{center}
    \begin{tabular}{c}
    \includegraphics[width=6.25cm]{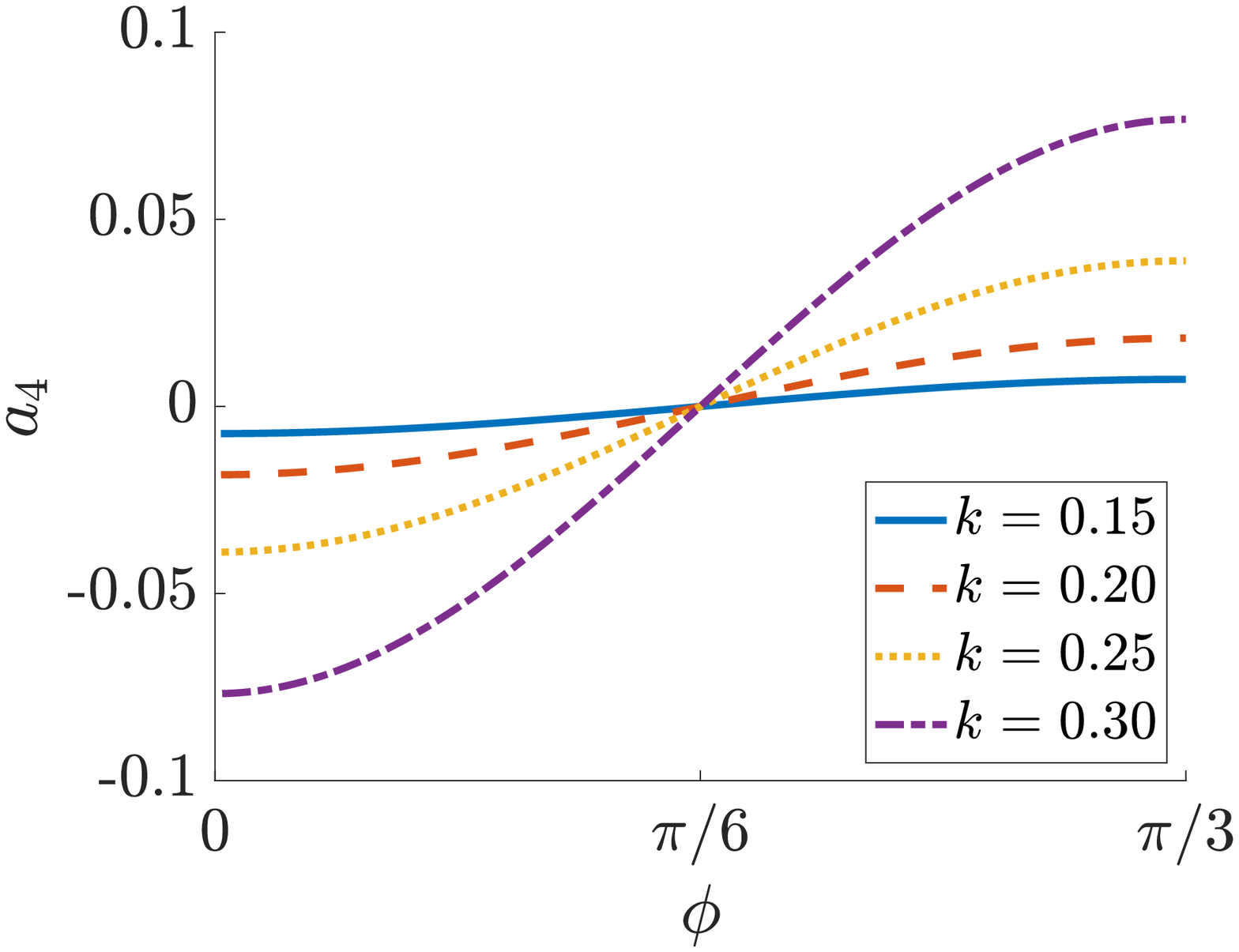} \\
    \includegraphics[width=6.25cm]{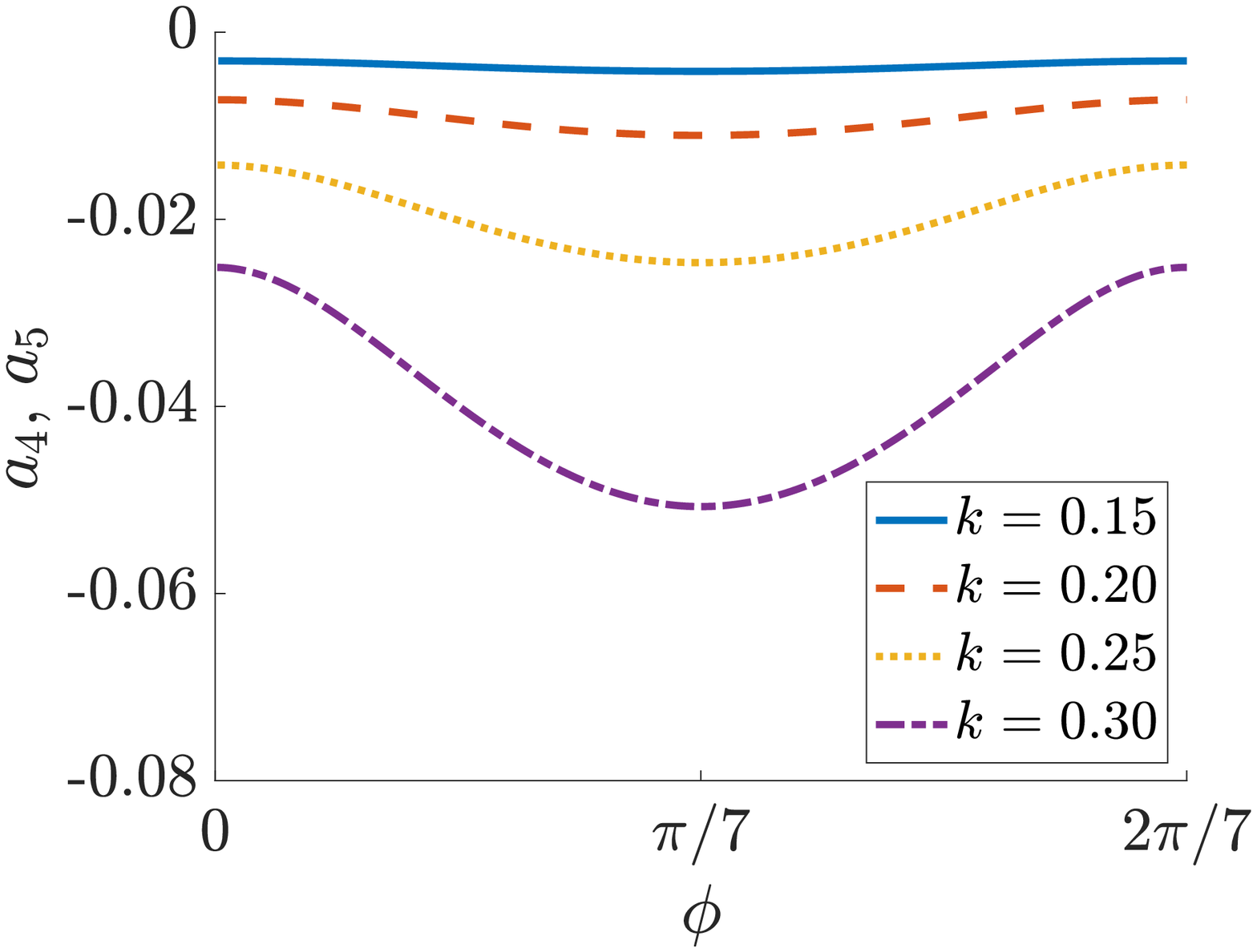} 
    \end{tabular}
    \end{center}
    \caption{Amplitude of node 4 (minimum intensity) versus $\phi$ for standing wave solution with $N=6$ (top), and amplitude of nodes 4 and 5 (minimum intensity) versus $\phi$ for standing wave solution with $N=7$ (bottom). $\omega = 1 \text{ mm}^{-1}$, coupling constants $k$ are $\text{mm}^{-1}$.}
    \label{fig:ABsuppression}
\end{figure}

The top panel of \cref{fig:ABsuppression} shows the amplitude of the node with minimum intensity (node 4 in \cref{fig:twist025}) versus the twist parameter $\phi$ for $N=6$. For all values of the coupling parameter $k$, the amplitude of this node is 0 when the twist parameter is given by $\phi = \pi/6$, which is an example of optical Aharonov-Bohm suppression. Since this a standing wave solution, this node will have 0 intensity for all $z$. This observation of a dark node opposite the node of maximum intensity agrees with the results of \cite{castro2016,Parto2017}. We show below in Section \ref{sec:Neven} that this occurs in general when $N$ is even and $\phi = \pi/N$. The bottom panel of \cref{fig:ABsuppression} shows the amplitude of the nodes with minimum intensity (nodes 4 and 5 in \cref{fig:twist025N7}) versus the twist parameter $\phi$ for $N=7$. Since the amplitudes of these nodes are never 0, optical Aharonov-Bohm suppression does not occur when $N$ is odd and there is a single excited node. (See Section \ref{sec:Nodd} below for a setting in which AB suppression does occur for $N$ odd).

\section{Optical Aharonov-Bohm suppression}\label{sec:ABsupp}

We now show that optical Aharonov-Bohm suppression occurs for standing wave solutions when the twist parameter is $\phi = \pi/N$. We consider the cases of $N$ even and $N$ odd separately, since the symmetry patterns are different. In both cases, we find that we can obtain a single dark node when $\phi = \pi/N$.

\subsection{\texorpdfstring{$N$}{N} even}\label{sec:Neven}

\begin{figure}
    \begin{center}
    \begin{tabular}{cc}
    \includegraphics[width=3.9cm]{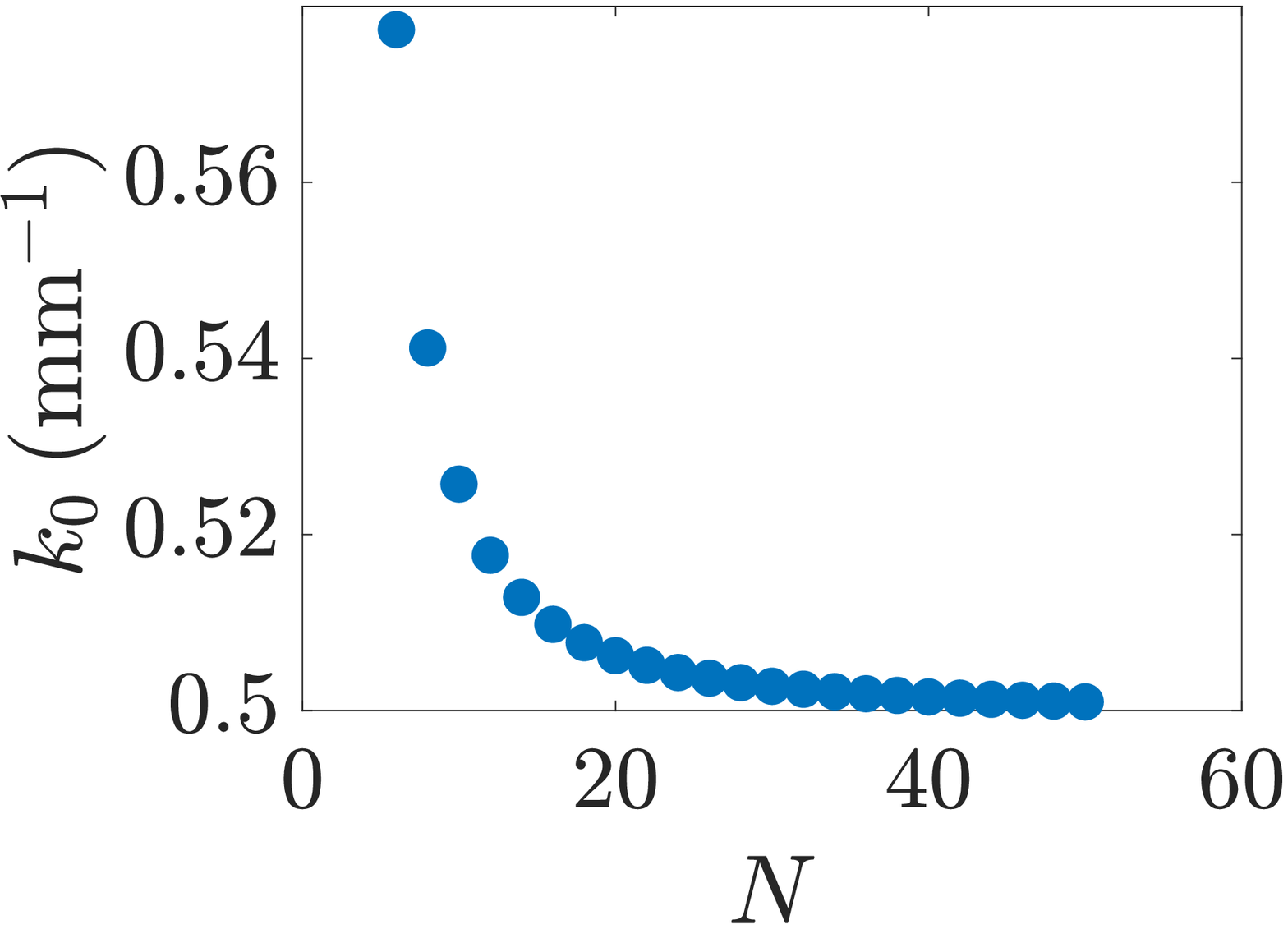} &
    \includegraphics[width=3.9cm]{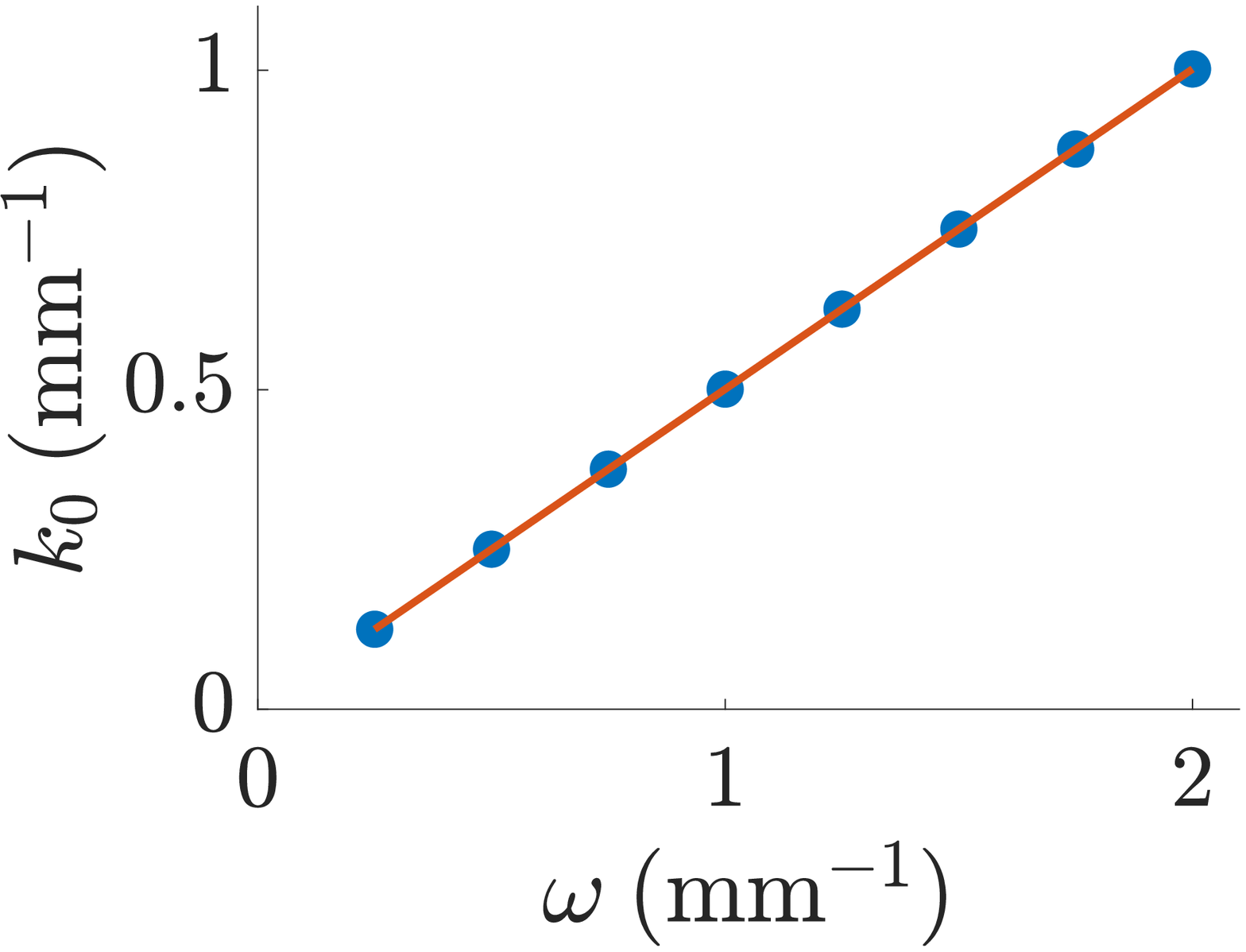}
    \end{tabular}
    \end{center}
    \caption{Left panel shows $k_0$ versus $N$ for dark node opposite peak intensity node for $N$ even, $\omega = 1 \text{ mm}^{-1}$. Right panel shows $k_0$ vs $\omega$ together with least squares linear regression line for $N = 50$.}
    \label{fig:k0plot}
\end{figure}

Taking $a_M = 0$, where $M = (N/2)+1$, we use the symmetries \cref{eq:symm} to reduce the system \cref{eq:twisteqreal} to
\begin{equation*}
\begin{aligned}
&2 k a_2 \cos(\theta_2 - \phi) + \omega a_1 - a_1^3 = 0 \\
&\begin{rcases}
k( a_{n+1} \cos(\theta_{n+1}-\theta_n-\phi) \\
\quad+a_{n-1} \cos(\theta_n - \theta_{n-1}-\phi)) \\
\quad+\omega a_n - a_n^3 = 0  \\
a_{n+1} \sin(\theta_{n+1}-\theta_n-\phi) \\
\quad- a_{n-1} \sin(\theta_n - \theta_{n-1}-\phi) = 0
\end{rcases} n = 2, \dots, M-1\\
&2 k a_{M-1} \cos(\theta_{M-1} + \phi) = 0 \\
& \theta_1 = \theta_M = 0.
\end{aligned}
\end{equation*}
It follows that $a_n = 0$ for all $n$ unless
\begin{equation*}
\begin{aligned}
&\cos(\theta_{M-1} + \phi) = 0 \\
&\sin(\theta_{n} - \theta_{n-1} - \phi) = 0 && \qquad n = 3, \dots, M-1 \\
&\sin(\theta_2 - \phi) = 0.
\end{aligned}
\end{equation*}
One solution to this is
\begin{equation*}
\begin{aligned}
&\theta_{M-1} + \phi = \pi/2 \\
&\theta_{n} - \theta_{n-1} - \phi = 0 && \qquad n = 3, \dots, M-1 \\
&\theta_2 - \phi = 0,
\end{aligned}
\end{equation*}
from which it follows that we can have a single dark node at site $M$ when $\phi = \pi/N$. If this is the case, the system of equations above reduces to the simpler system of $M-1$ equations
\begin{equation}\label{eq:twisteqevenhole}
\begin{aligned}
&2 k a_2 + \omega a_1 - a_1^3 = 0 \\
&k\left( a_{n+1} + a_{n-1} \right) + \omega a_n - a_n^3 = 0 \quad n = 2, \dots, M-2 \\
&k a_{M-2} + \omega a_{M-1} - a_{M-1}^3 = 0.
\end{aligned}
\end{equation}
This system is of the form $F(a,k) = 0$, where $a = (a_1, \dots, a_{M-1})$. $F(\tilde{a}, 0) = 0$, where $\tilde{a} = (\sqrt{\omega}, 0, \dots, 0)$. Since $D_F(\tilde{a}, 0) = \diag(-2\omega,\omega, \dots, \omega)$, which is invertible for $\omega \neq 0$, it follows from the implicit function theorem that there exists $k_0 > 0$ such the system \cref{eq:twisteqevenhole} has a unique solution for all $k$ with $|k| < k_0$. The critical value $k_0$ can be computed numerically by parameter continuation with AUTO, and will depend on both $N$ and $\omega$ (\cref{fig:k0plot}). These computations suggest that $k_0$ approaches $\omega/2$ as $N$ becomes large.

\begin{figure}
\begin{center}
\begin{tabular}{c}
\includegraphics[width=6.25cm]{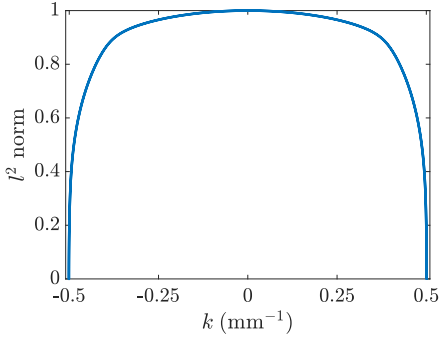} 
\end{tabular}
\end{center}
\caption{$\ell^2$ norm of solution vs $k$ for $N = 50$ with dark node opposite peak intensity node. $\omega = 1 \text{ mm}^{-1}$.}
\label{fig:evenbif}
\end{figure}

\cref{fig:evenbif} plots the $\ell^2$ norm of the solution to \cref{eq:twisteqevenhole}
\begin{equation}
\| a \|_{\ell^2} = \left( \sum_{j=1}^{N/2} |a_j|^2 \right)^{1/2}
\end{equation}
versus the coupling parameter $k$. The critical value $k_0$ is the point at which the bifurcation curve touches the horizontal axis. As $k$ approaches $k_0$ in the parameter continuation, the $\ell^2$ norm of the solution approaches 0, thus the solution approaches the zero solution. Although it is possible that there are standing wave solutions for $|k| > k_0$, they cannot be reached by parameter continuation from this branch of solutions. At $k = 0$ (the AC limit), there is only one excited node with intensity $\sqrt{\omega}$, thus the $\ell^2$ norm of that solution is $\sqrt{\omega}$ (in \cref{fig:evenbif}, $\omega = 1$, thus the $\ell^2$ norm of the solution is 1 when $k = 0$).

\begin{figure}
\begin{center}
\begin{tabular}{c}
\includegraphics[width=8cm]{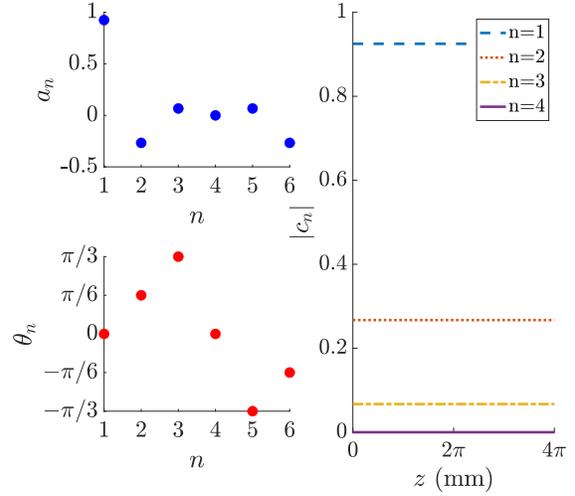}
\end{tabular}
\end{center}
\caption{Standing wave solution for $N = 6$ and $\phi = \pi/6$. Left panel shows amplitudes $a_n$ and phases $\phi_n$ for solution at each node. Right panel is intensity of solution $|c_n|$ versus $z$ for nodes 1-4, which is constant in $z$. Evolution in $z$ computed using fourth order Runge-Kutta method. Node 1 has maximum amplitude, and node 4 is a dark node. $\omega = 1 \text{ mm}^{-1}$, $k = 0.25 \text{ mm}^{-1}$.}
\label{fig:evenhole6}
\end{figure}

Once \cref{eq:twisteqevenhole} has been solved numerically, the full solution to \cref{eq:twisteqreal} is given by
\begin{align*}
&a_M = 0 \\
&a_{M+j} = a_{M-j} && \qquad j = 1, \dots, M-2 \\
&\theta_0 = 0 \\
&\theta_j = (j-1)\phi && \qquad  j = 2, \dots, M-1 \\
&\theta_M = 0 \\
&\theta_{M+j} = -\theta_{M-j} && \qquad j = 1, \dots, M-2.
\end{align*}
\cref{fig:evenhole6} shows this solution for $N=6$ and $k=\pi/6$. The amplitudes $a_n$ and phases $\theta_n$ are qualitatively similar to the case when $\phi = 0.25$ (\cref{fig:twist025}); however, when $\phi = \pi/6$, the intensity of node 4 is equal to 0, whereas for other values of $\phi$, the intensity of node 4 is small, but nonzero (\cref{fig:ABsuppression}).

In \cite{castro2016}, a perturbation method is used for the $N=6$ case to show that if the peak intensity is contained in node 1, the opposite node (node 4) has an intensity of 0, to leading order, when $\phi=\pi/6$. Our analysis confirms the result of these asymptotics, but is much more rigorous in that it demonstrates that for all $N$ even, when the twist is given by $\phi = \pi/N$, a standing wave solution exists for which the peak intensity is contained in a single node, and the opposite node in the ring has intensity identically equal to 0 for all $z$.

\subsection{\texorpdfstring{$N$}{N} odd}\label{sec:Nodd}

\begin{figure}
\begin{center}
\begin{tabular}{c}
\includegraphics[width=8cm]{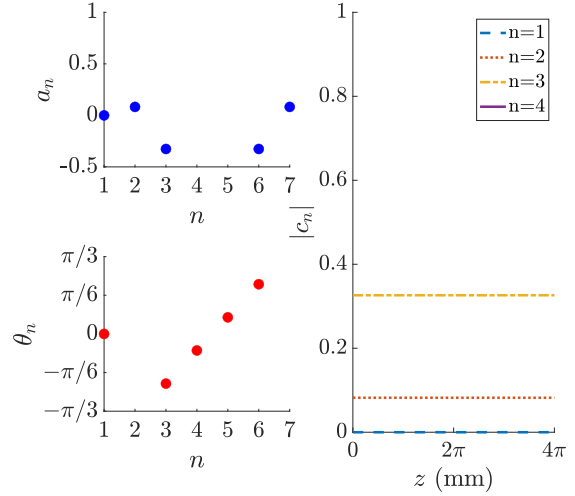}
\end{tabular}
\end{center}
\caption{Standing wave solution for $N = 7$ and $\phi = \pi/7$. Left panel shows amplitudes $a_n$ and phases $\phi_n$ for solution at each node. Right panel is intensity of solution $|c_n|$ versus $z$ for nodes 1-4, which is constant in $z$. Evolution in $z$ computed using fourth order Runge-Kutta method. Nodes 4 and 5 have equal and maximum amplitude, and node 1 is a dark node. $\omega = 1 \text{ mm}^{-1}$, $k = 0.25 \text{ mm}^{-1}$.}
\label{fig:oddhole7}
\end{figure}

When $N$ is odd and the peak intensity is contained in a single node, we cannot obtain dark nodes for any value of the twist parameter $k$ (\cref{fig:ABsuppression}). We can, however, obtain a dark node when $N$ is odd if we start with two adjacent bright nodes. For simplicity, we take node 1 to be the dark node; in this case, the dark node will be opposite a pair of bright nodes at $a_M$ and $a_{M+1}$ with the same amplitude, where $M = (N+1)/2$. Using the symmetries \cref{eq:symm}, when $a_1 = 0$, the system \cref{eq:twisteqreal} reduces to 
\begin{equation*}
\begin{aligned}
&2 k a_2 \cos(\theta_2 - \phi) = 0 \\
&k a_3 \cos(\theta_3-\theta_2-\phi) + \omega a_2 - a_2^3 = 0 \\
&a_3 \sin(\theta_3-\theta_2-\phi) = 0 \\
&\begin{rcases}
k( a_{n+1} \cos(\theta_{n+1}-\theta_n-\phi) \\
\quad+ a_{n-1} \cos(\theta_n - \theta_{n-1}-\phi)) \\
\quad+ \omega a_n - a_n^3 = 0 \\
a_{n+1} \sin(\theta_{n+1}-\theta_n-\phi) \\
\quad- a_{n-1} \sin(\theta_n - \theta_{n-1}-\phi) 
\end{rcases} n = 3, \dots, M-1 \\
&k ( a_M \cos(-2 \theta_M - \phi) + a_{M-1} \cos(\theta_M - \theta_{M-1} - \phi)) \\
&\quad+ \omega a_M - a_M^3 = 0 \\
& a_M \sin(-2 \theta_M - \phi) - a_{M-1} \sin(\theta_M - \theta_{M-1} - \phi) = 0.
\end{aligned}
\end{equation*}
It follows that $a_n = 0$ for all $n$ unless
\begin{equation*}
\begin{aligned}
&\cos(\theta_2 - \phi) = 0 \\
&\sin(\theta_{n} - \theta_{n-1} - \phi) = 0 && \qquad n = 3, \dots, M-1 \\
&\sin(2 \theta_M + \phi) = 0.
\end{aligned}
\end{equation*}
One solution to this is
\begin{equation}\label{eq:odddarknodecond1}
\begin{aligned}
&\theta_2 - \phi = -\pi/2 \\
&\theta_{n} - \theta_{n-1} - \phi = 0 && \qquad n = 3, \dots, M-1 \\
&2 \theta_M + \phi = 0,
\end{aligned}
\end{equation}
from which it follows that we can have a single dark node at $a_1$ when $\phi = \pi/N$. This condition for a single dark node is the same as when $N$ is even. For this case, the system of equations above reduces to the simpler system of equations
\begin{equation}\label{eq:twisteqoddhole}
\begin{aligned}
& k a_3 + \omega a_2 - a_2^3 = 0\\
&k( a_{n+1} + a_{n-1} ) + \omega a_n - a_n^3 = 0 \quad n = 3, \dots, M-1 \\
&k ( a_M + a_{M-1} ) + \omega a_M - a_M^3 = 0.
\end{aligned}
\end{equation}
This system of equations is again of the form $F(a,k) = 0$, where $a = (a_2, \dots, a_M)$. $F(\tilde{a}, 0) = 0$, where $\tilde{a} = (0, \dots, 0, \sqrt{-\omega/d}, 0)$. Since $D_F(\tilde{a}, 0) = \diag(\omega, \dots, \omega, -2\omega)$, which is invertible for $\omega \neq 0$, it follows from the implicit function theorem that there exists $k_0 > 0$ such the system \cref{eq:twisteqoddhole} has a unique solution for all $k$ with $|k| < k_0$. As in the case for $N$ even, the critical value $k_0$, as well as its dependency on $N$ and $\omega$, can be computed numerically. Once \cref{eq:twisteqoddhole} has been solved numerically, we obtain the full solution to \cref{eq:twisteqreal} using
\begin{align*}
&a_1 = 0 \\
&a_{M+j} = a_{M-j+1} && \qquad j = 1, \dots, M-1 \\
&\theta_0 = 0 \\
&\theta_j = (j-1)\phi - \pi/2 && \qquad j = 2, \dots, M \\
&\theta_{M+j} = -\theta_{M-j+1} && \qquad j = 1, \dots, M-1
\end{align*}
\cref{fig:oddhole7} shows this solution for $N=7$. The peak intensity in this solution is contained in two adjacent nodes, and there is a single dark node opposite this pair, which is qualiatively different from the solution in \cref{fig:twist025N7}.

\section{Stability}\label{sec:stability}

\begin{figure}
\begin{center}
\begin{tabular}{cc}
\includegraphics[width=4cm]{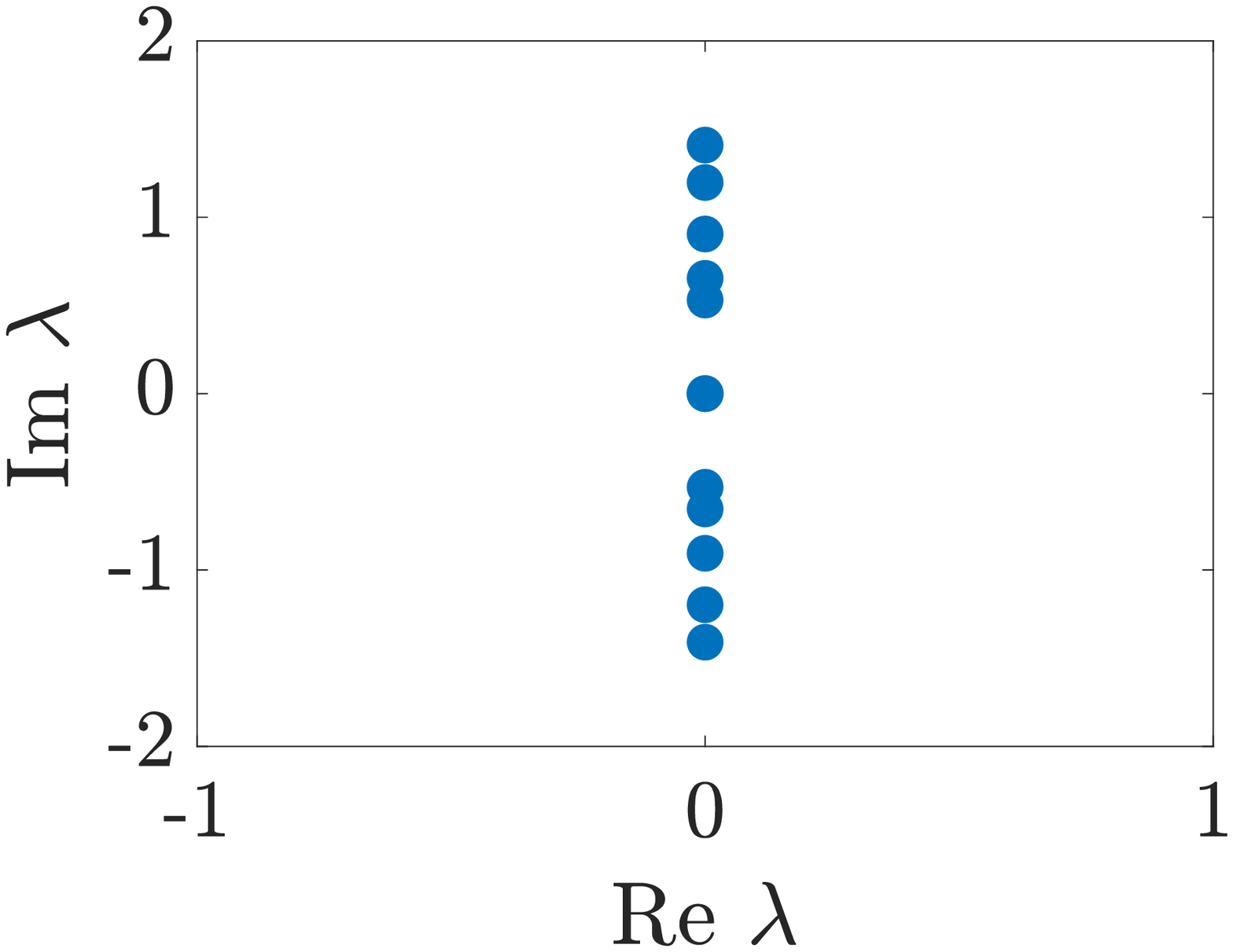}
\includegraphics[width=4cm]{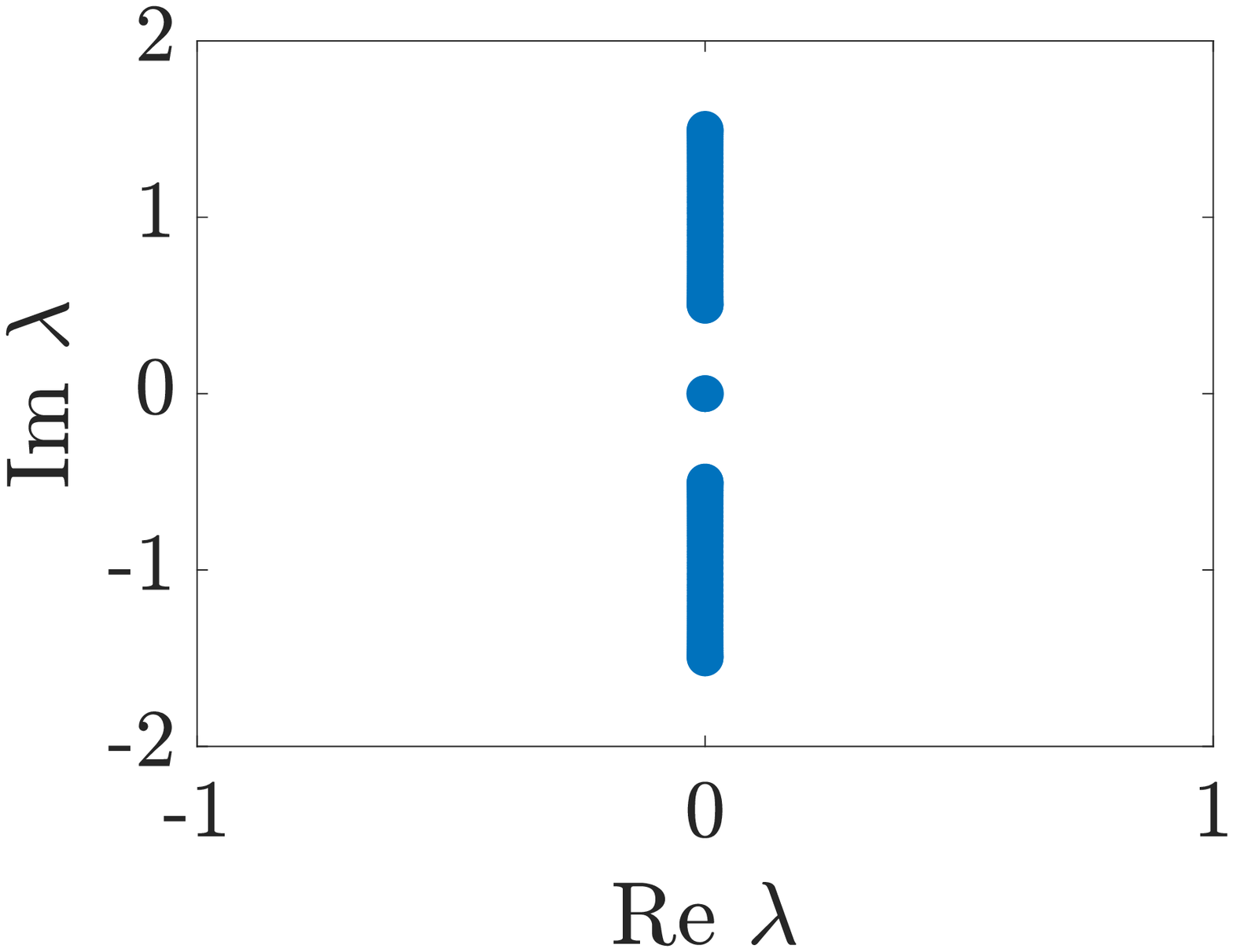}
\end{tabular}
\end{center}
\caption{Spectrum of linearization of \cref{eq:twist1} about solution for even $N$ with a single dark node opposite a single bright node. $N=6$ (left panel) and $N=50$ (right panel). $k=0.25\text{ mm}^{-1}$, $\omega = 1\text{ mm}^{-1}$, $\phi = \pi/N$.}
\label{fig:evenholespec}
\end{figure}

We now look at the stability of the standing wave solutions we constructed in the previous section. As a first step in stability analysis, the linearization of equation \cref{eq:twist1} about a standing wave solution $c_n = a_n e^{i (\omega z + \theta_n) } = (v_n + i w_n)e^{i\omega z}$ is the $2N \times 2N$ block matrix
\begin{equation}\label{eq:linearization}
\begin{aligned}
A(&c_n) =
k \begin{pmatrix}S & C \\ -C & S \end{pmatrix}
+ \omega\begin{pmatrix}0 & I \\ -I & 0 \end{pmatrix} \\
&-\begin{pmatrix} \diag(2v_n w_n) & \diag(v_n^2 + 3 w_n^2) \\
-\diag(3 v_n^2 + w_n^2) & -\diag(2v_n w_n) \end{pmatrix}
\end{aligned}
\end{equation}
where each block is a $N\times N$ matrix, $C$ is the periodic banded matrix with $\cos \phi$ on the first upper and lower diagonals, and $S$ is the periodic banded matrix with $\sin \phi$ on the first lower diagonal and $-\sin \phi$ on the first upper diagonal, i.e.
\begin{align*}
C &= \begin{pmatrix}
0 & \cos \phi & & \dots & \cos \phi \\
\cos \phi & 0 & \cos \phi & & & \\
& & \ddots & \ddots &  & \\
\cos \phi & & \dots & \cos \phi & 0
\end{pmatrix} \\ 
S &= \begin{pmatrix}
0 & -\sin \phi & & \dots & \sin \phi \\
\sin \phi & 0 & -\sin \phi & & & \\
& & \ddots & \ddots &  & \\
-\sin \phi & & \dots & \sin \phi & 0
\end{pmatrix}.
\end{align*}
Since \cref{eq:linearization} is a finite dimensional matrix, the spectrum is purely point spectrum. Due to the gauge invariance, there is an eigenvalue at 0 with algebraic multiplicity 2 and geometric multiplicity 1. Following the analysis in \cite[Section 2.1.1.1]{Kevrekidis2009}, there are plane wave eigenfunctions which are, to leading order, of the form $e^{\pm( i q n + \lambda z)}$, where $q$ is the discrete wavenumber, and satisfy the dispersion relation
\begin{equation}\label{eq:dispersion}
\lambda = \pm i \left( \omega + 2 k \cos(q + \phi) \right).
\end{equation}
The corresponding eigenvalues $\lambda$ are thus purely imaginary and are contained in the bounded intervals $\pm i[\omega - 2 k, \omega + 2 k]$. As $N$ increases, these eigenvalues fill out this interval. For $|k| < k_0 = \omega/2$, these eigenvalues do not interact with the kernel eigenvalues. \cref{fig:evenholespec} illustrates these results numerically for $\omega = 1$ and $k = 0.25$ for the case of $N$ even and $\phi = \pi/N$, i.e. a single dark node opposite a single bright node. Similar results are obtained for other values of $\omega$ and $k$ in which there is a single bright node as well as the solutions from Section \ref{sec:Nodd} with odd $N$ and a single dark node opposite a pair of bright nodes.

\begin{figure}
\begin{center}
\begin{tabular}{cc}
\includegraphics[width=4cm]{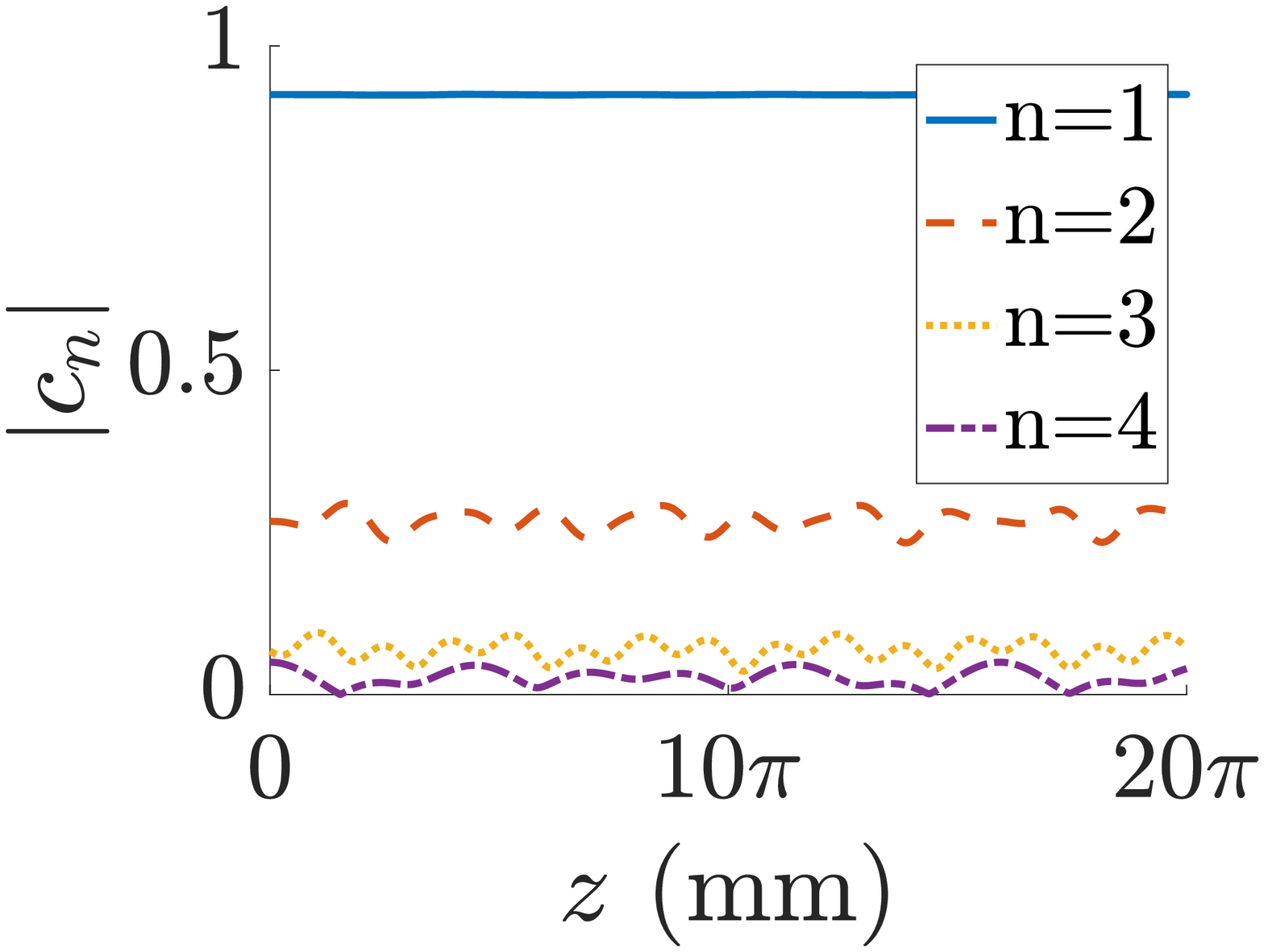} & 
\includegraphics[width=4cm]{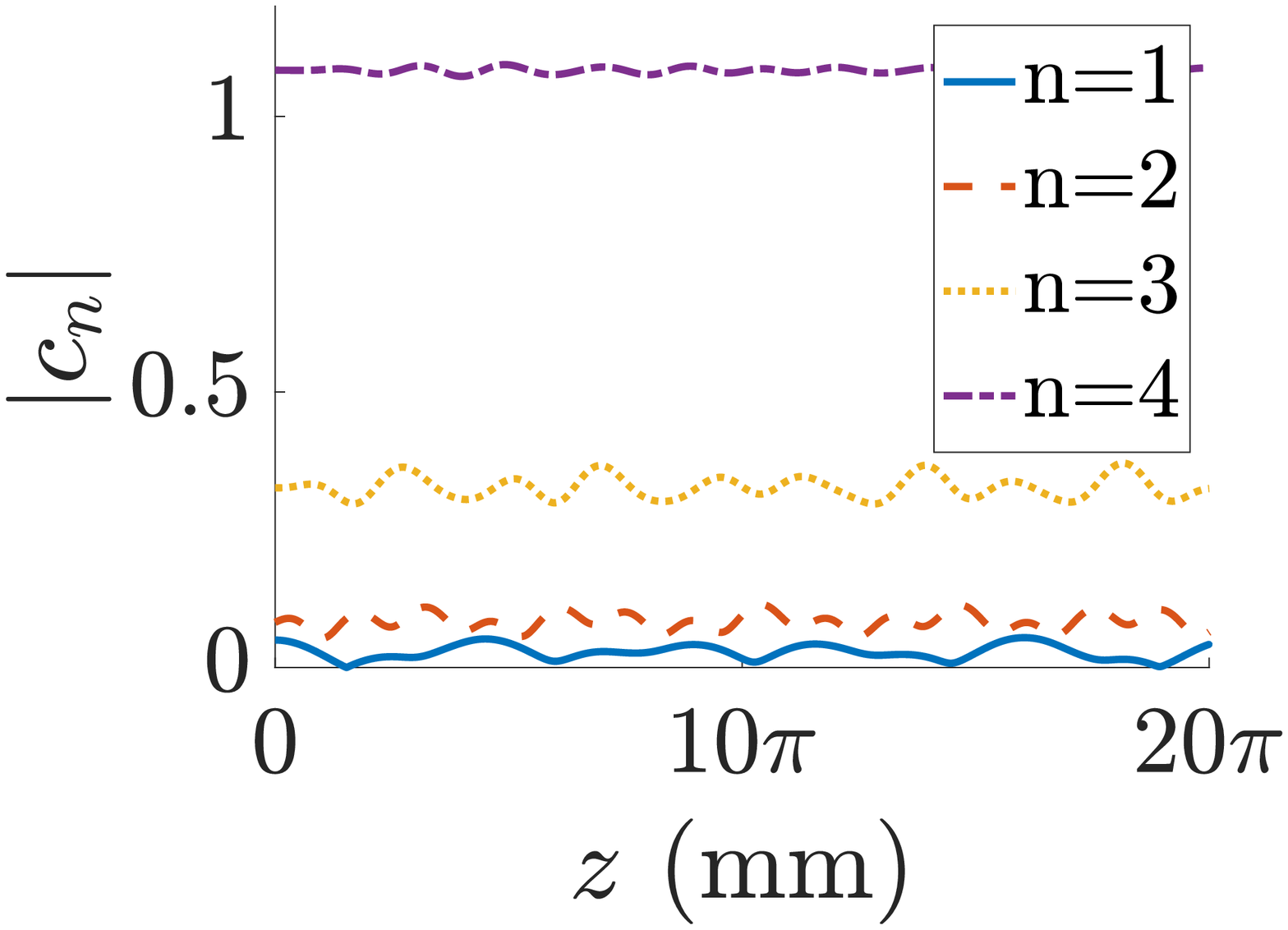}
\end{tabular}
\end{center}
\caption{Amplitude $|c_n|$ for first four nodes versus $z$ for solution with $N=6$, $\phi = \pi/6$ (left panel) and $N=7$, $\phi = \pi/7$ (right panel). Evolution performed using a fourth order Runge-Kutta scheme, $k=0.25 \text{ mm}^{-1}$.}
\label{fig:evenhole6perturbed}
\end{figure}

Since the spectrum of these solutions is purely imaginary, we expect that they will be neutrally stable, i.e. any small perturbation will remain close to the the unperturbed standing wave for all $z$, but will exhibit oscillatory behavior. \cref{fig:evenhole6perturbed} shows the results of numerical evolution in $z$ for a small perturbation of the standing wave solution when $N=6$ and $N=7$. For the initial condition of the perturbed solution, a small quantity (0.05) was added to the amplitude of dark node. (This initial condition was chosen for simplicity. Any initial condition which is close to the unperturbed solution in amplitude and phase produces results which are qualitatively the same). Numerical results show small amplitude oscillations about the unperturbed solutions, but no growth, which provides numerical evidence for neutral stability. The amplitude of the oscillations depends on the magnitude of the initial perturbation. (Compare these evolution plots to the right panels of \cref{fig:evenhole6} and \cref{fig:oddhole7}, noting that the evolution in $z$ in \cref{fig:evenhole6perturbed} is over a much greater length). Similar results are obtained for other values of $N$, $\omega$, and $k$.

\begin{figure}
\begin{center}
\begin{tabular}{cc}
\includegraphics[width=4cm]{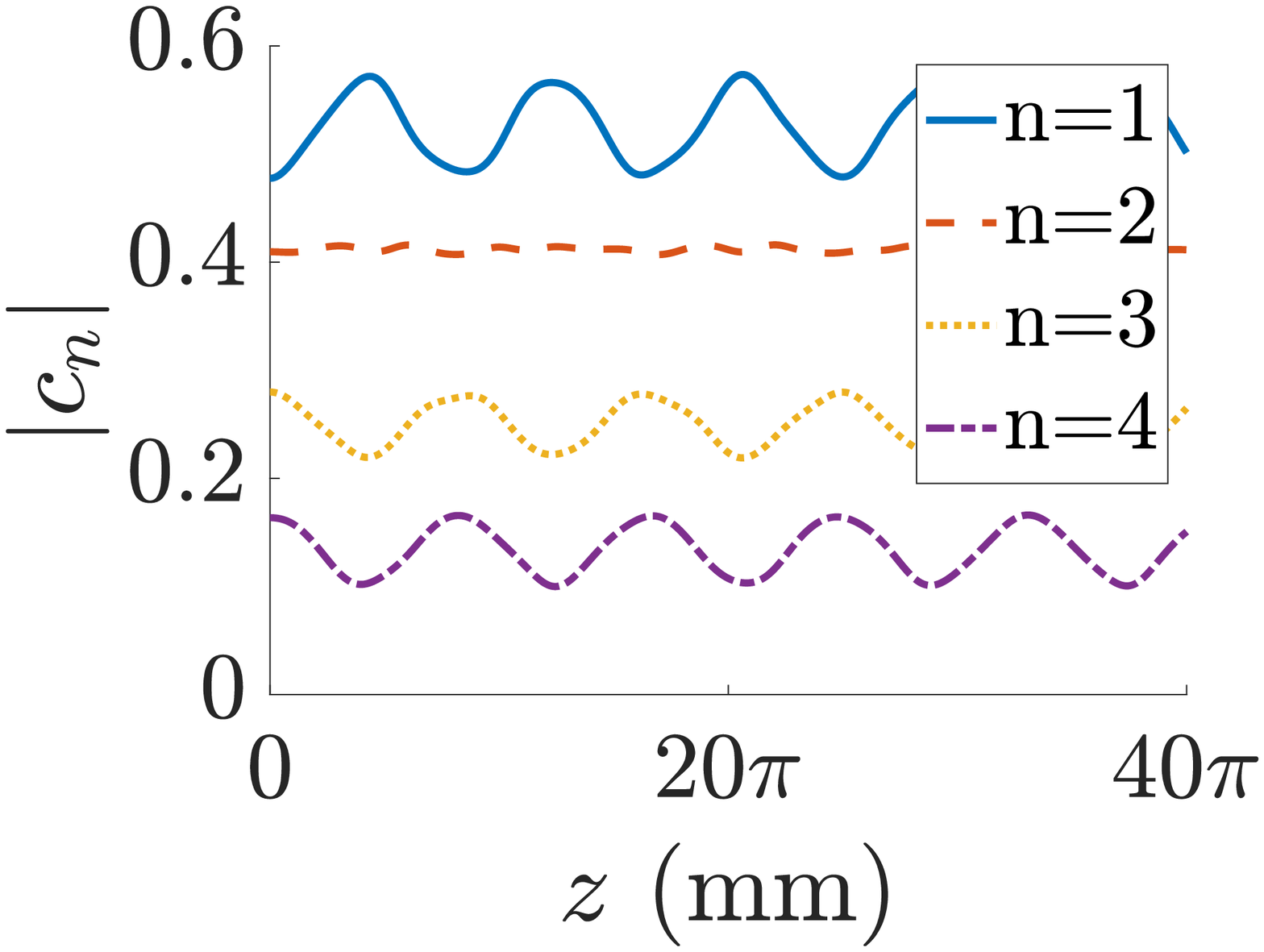} & 
\includegraphics[width=4cm]{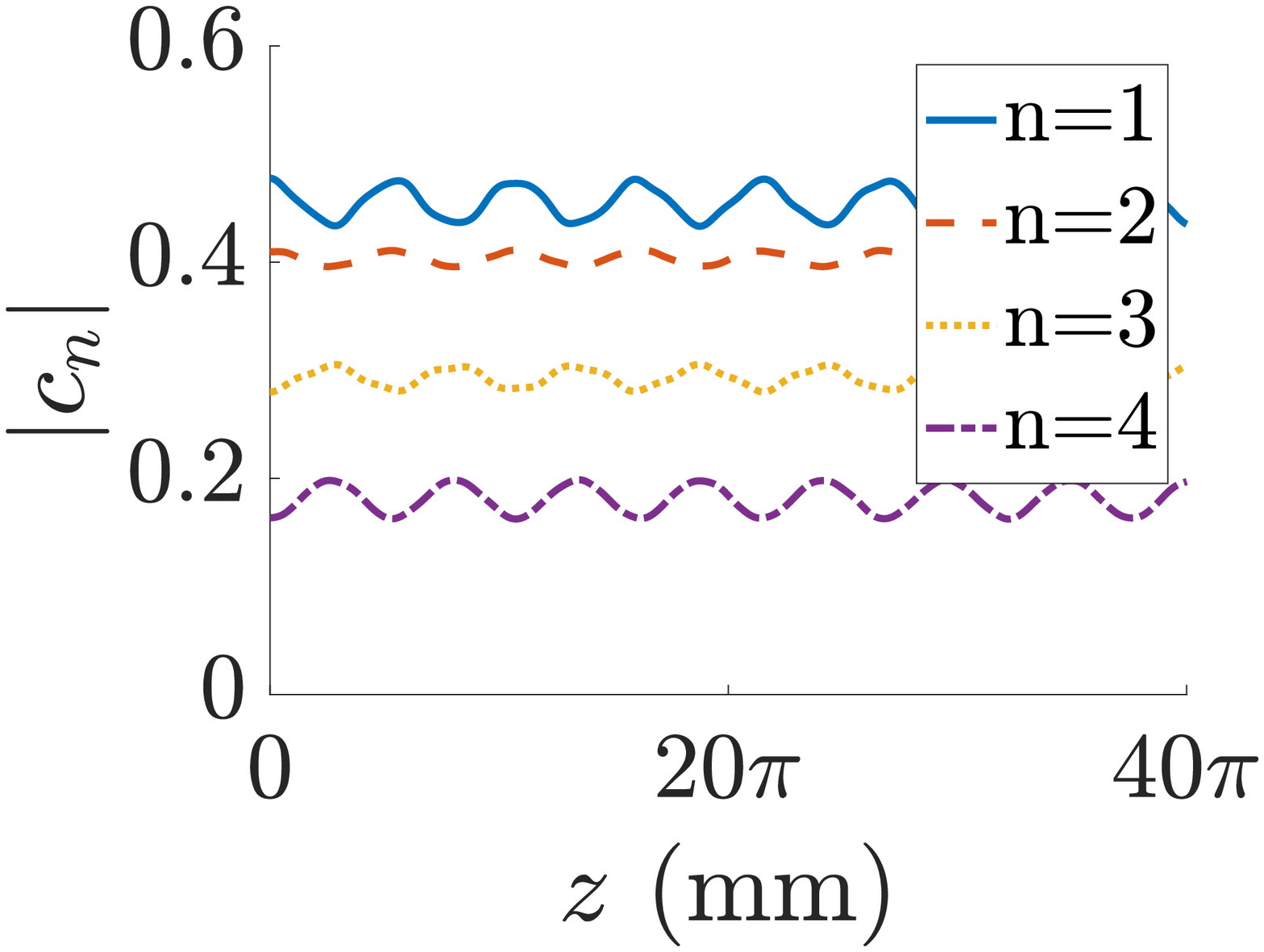}
\end{tabular}
\end{center}
\caption{Amplitude $|c_n|$ for first four nodes versus $z$ for solution with $N=10$ and $\phi = \pi/10$. Initial condition is solution to \cref{eq:twisteqevenhole} with $k = 0.45\text{ mm}^{-1}$. Evolution performed with $k = 0.35\text{ mm}^{-1}$ (left) and $k = 0.55\text{ mm}^{-1}$ (right) using a fourth order Runge-Kutta scheme.}
\label{fig:evenholekperturbed}
\end{figure}

In addition, we can start with a neutrally stable standing wave solution and perturb the system by a small change in $k$ or $\phi$. \cref{fig:evenholekperturbed} shows the results of perturbations in $k$. In particular, note that in the right panel of \cref{fig:evenholekperturbed}, the system is evolved using a value of the coupling parameter $k$ which is greater than $k_0$, where $k_0$ is defined in Section \ref{sec:Neven}. In both cases, the solutions show oscillations, indicating this to be robust dynamics. The simulation suggests the period of oscillations has a strong dependence on $k$. Additional evolution results can be found in \cite{castro2016}. In particular, see \cite[Figure 4]{castro2016} for evolution results when the fiber is initially excited at a single site. 

\section{Asymmetric coupling}\label{sec:variants}

As an additional variant, if the strength of the nearest-neighbor coupling is allowed to differ between pairs of nodes, equation \cref{eq:twist1} becomes
\begin{equation}\label{eq:twistk}
i \partial_z c_n = k_{n+1} e^{-i\phi}c_{n+1} + k_{n-1} e^{i\phi}c_{n-1} \pm |c_n|^2 c_n,
\end{equation}
where there is a different coupling parameter $k_n$ for each pair of nodes. As in the symmetric case, we will  only consider the defocusing (minus) nonlinearity. An asymmetric configuration can be realized by either having a variation in the separation between waveguides or having a variation of the core radius, although in the latter case, variations of the propagation constant must be accounted by adding a term of the form $k_n c_n$ to \cref{eq:twistk}. Even in the idealized case of identical separation, small variations could appear as a consequence of imperfections in the fiber bundle construction, in which case the parameters $k_n$ would be close, but not identical. 
\begin{figure}
\begin{center}
\begin{tabular}{c}
\includegraphics[width=8cm]{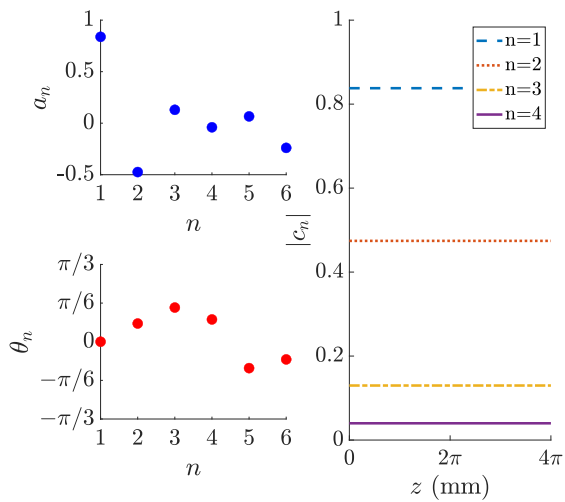}
\end{tabular}
\end{center}
\caption{Asymmetric standing wave solution to \cref{eq:twistk} for $N = 6$, $k_1 = 0.4\text{ mm}^{-1}$, and $k_n = 0.25\text{ mm}^{-1}$ for all other $n$. Left panel shows amplitudes $a_n$ and phases $\phi_n$ for solution at each node. Right panel is intensity of solution $|c_n|$ versus $z$ for nodes 1-4, which is constant in $z$. Evolution in $z$ computed using fourth order Runge-Kutta method. $\phi = 0.25$, $\omega = 1\text{ mm}^{-1}$.}
\label{fig:even6assym}
\end{figure}
\noi This allows for asymmetric solutions, as shown in \cref{fig:even6assym}. When compared to the symmetric solution for uniform $k$ in \cref{fig:twist025} (which has the same set of parameters except for the coupling parameter $k_1$), the phases and amplitudes are similar in magnitude, but the symmetry relations \cref{eq:symm} have been broken. The asymmetric solution in \cref{fig:even6assym} is neutrally stable, since its spectrum is imaginary, and small perturbations result in oscillatory behavior about the unperturbed solution (\cref{fig:assymstab}, compare to the right panel of \cref{fig:even6assym}). Although a rigorous analysis of these asymmetric solutions is beyond the scope of this work, the results of these numerical simulations suggest that it is likely that small differences in the coupling coefficients $k_n$ do not affect stability, which would imply that small imperfections in the physical model would not result in loss of stability.

\begin{figure}
    \begin{center}
    \begin{tabular}{cc}
    \includegraphics[width=4cm]{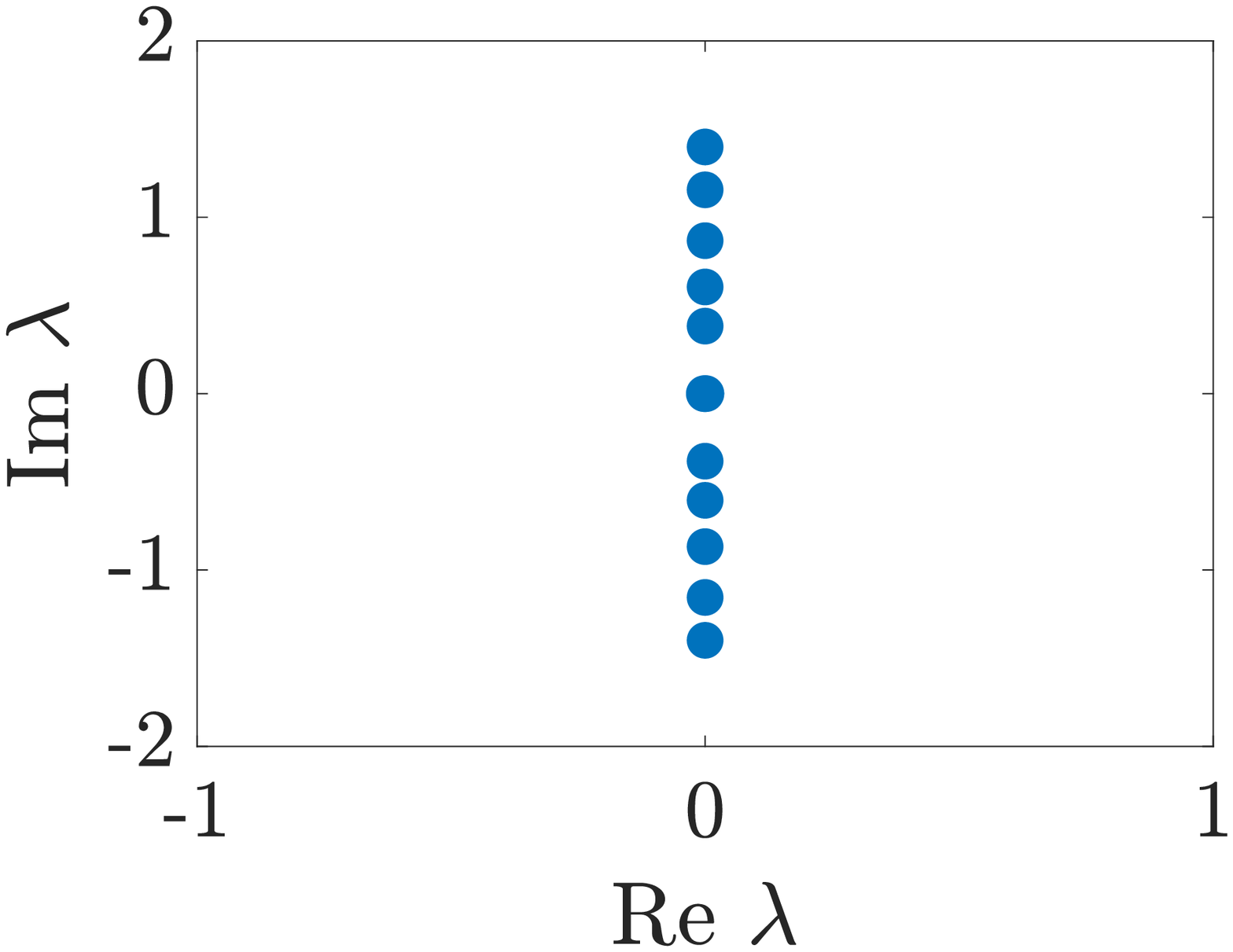}
    \includegraphics[width=4cm]{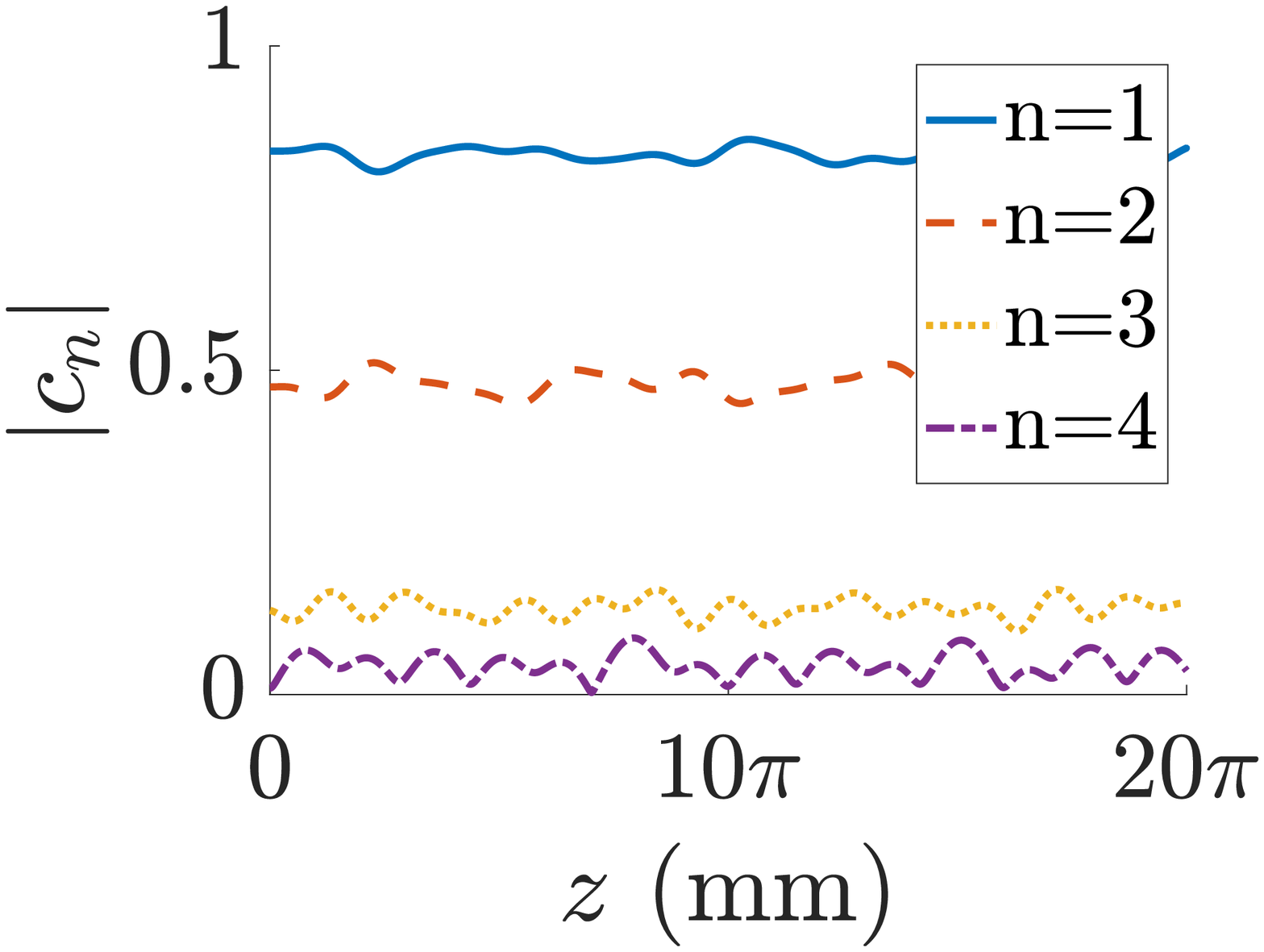}
    \end{tabular}
    \end{center}
    \caption{Spectrum of linearization of \cref{eq:twistk} about asymmetric standing wave solution from \cref{fig:even6assym} (left panel). Amplitude $|c_n|$ for first four nodes versus $z$ for evolution of perturbation of this solution using fourth order Runge-Kutta scheme (right panel).}
    \label{fig:assymstab}
    \end{figure}

\section{Multi-pulses}

\begin{figure}
\begin{center}
\begin{tabular}{cc}
\includegraphics[width=4cm]{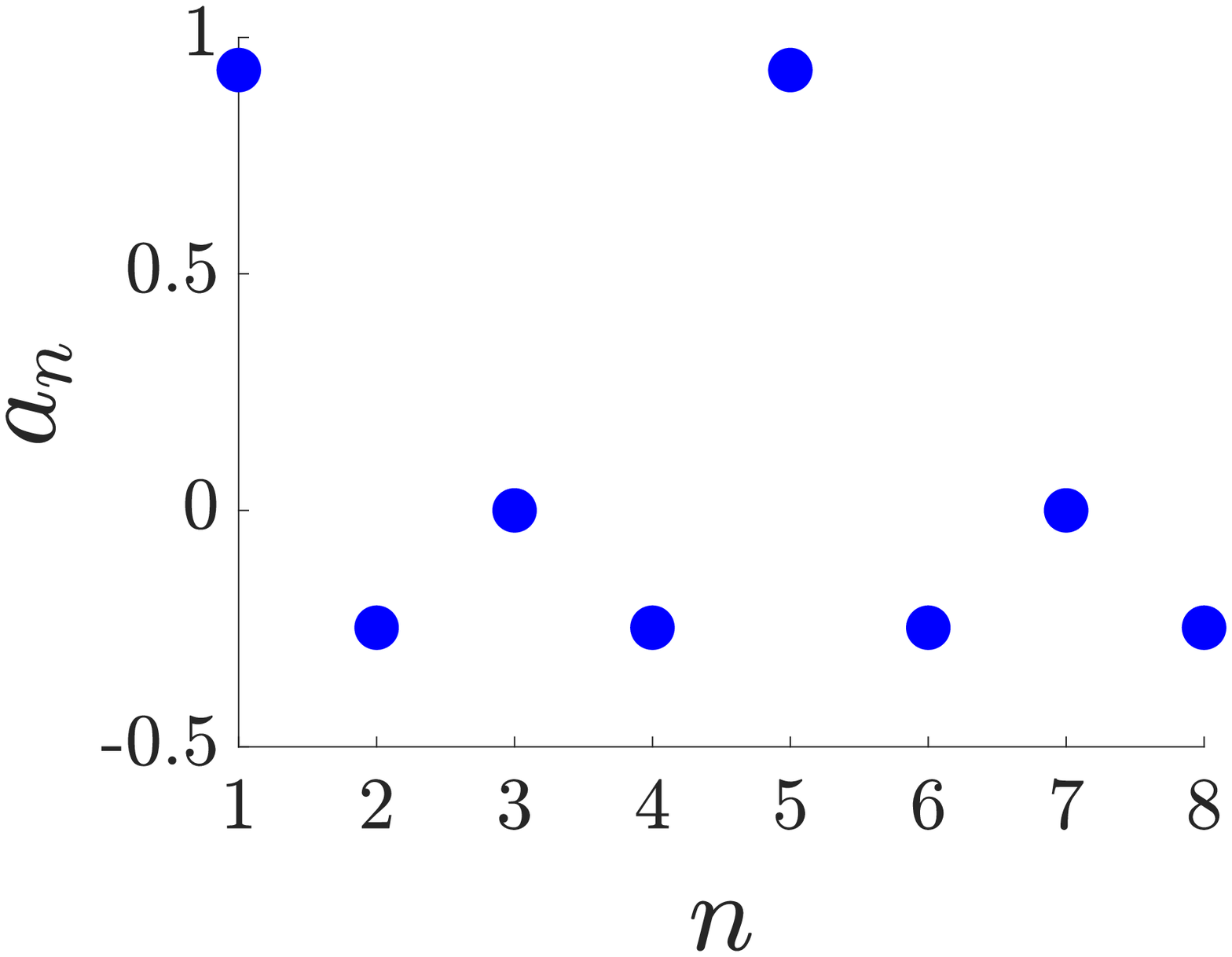} &
\includegraphics[width=4cm]{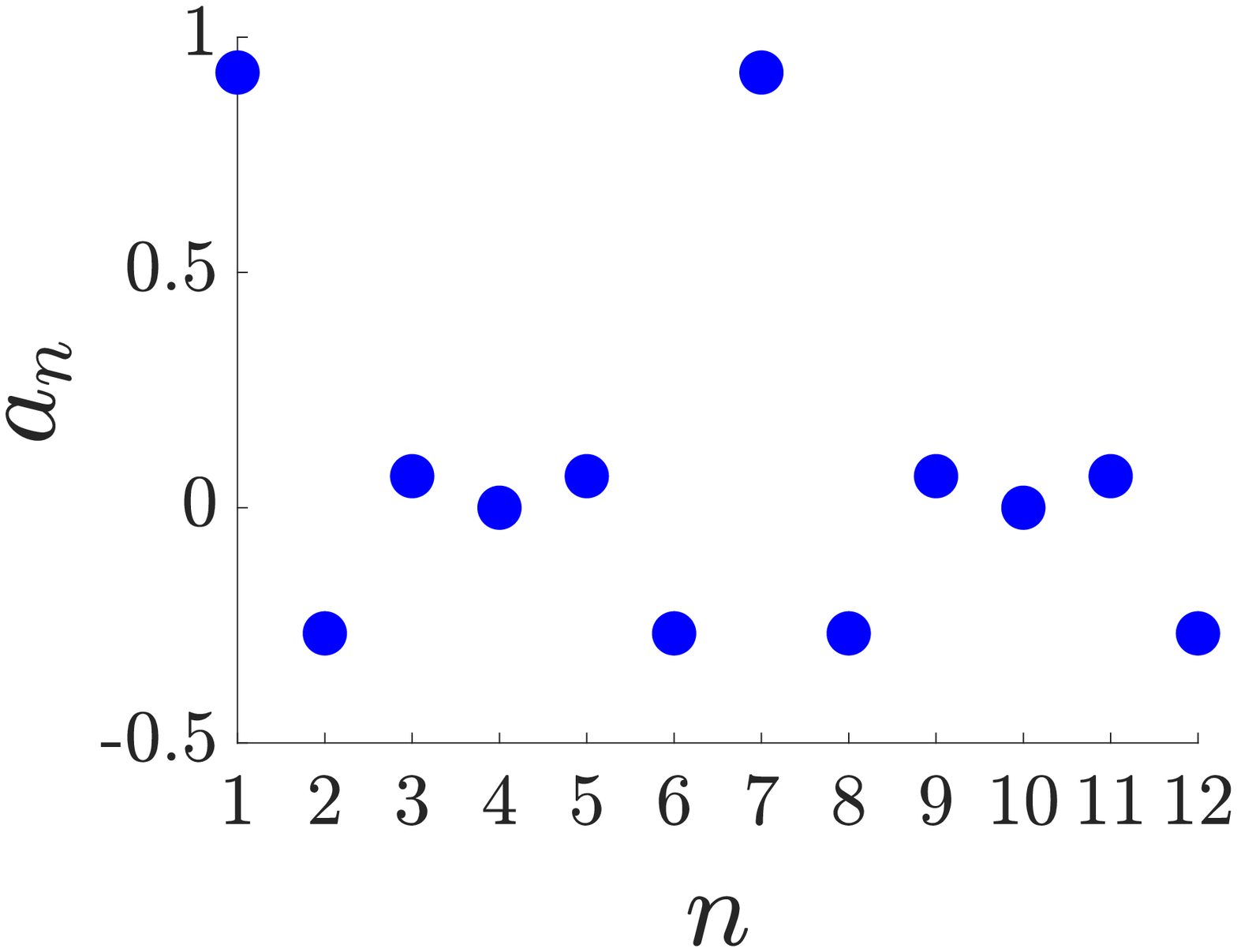}
\end{tabular}
\end{center}
\caption{Amplitudes $a_n$ for double pulse solutions with two bright nodes in opposite positions of the ring. $N=8$, $\phi=\pi/4$ (left panel) and $N=12$, $\phi=\pi/6$ (right panel). $\omega = 1\text{ mm}^{-1}$, $k=0.25\text{ mm}^{-1}$.}
\label{fig:dp}
\end{figure}

Another broad class of solutions is multi-pulses, which are solutions in which the energy is concentrated at multiple nodes which are well separated in the ring (see \cref{fig:dp} for two examples where the two nodes with peak intensity occupy opposite positions). In contrast with the solutions in \cref{fig:oddhole7}, where the intensity is concentrated at two adjacent nodes, the energy in a multi-pulse is concentrated at sites which are far apart. The solutions with two adjacent excited sites behave like a single soliton (see \cite{Kevrekidis2009} for a discussion of on-site and intersite solitons in the discrete NLS equation), whereas multi-pulses behave like a collection of solitons which can interact with their neighbors on either side \cite{Parker2020}.

Multi-pulses can be generated by parameter continuation from the AC limit, similar to what was done in \cref{sec:standingwave}. Although a systematic study of the existence and stability of multi-pulses is beyond the scope of this paper (see, for example, \cite{Parker2020} for results on multi-pulses in the discrete NLS equation), we present one example of a symmetric double pulse solution for even $N$ in which the two excited sites are opposite each other in the ring. If $N$ is a multiple of 4 and $\phi = 2\pi/N$, there is a pair of dark nodes halfway between the two bright nodes (in both directions), as can be seen in \cref{fig:dp}. In fact, these particular solutions are exactly two copies of the solutions from Section \ref{sec:Neven} spliced together.

\begin{figure}
\begin{center}
\begin{tabular}{cc}
\includegraphics[width=4cm]{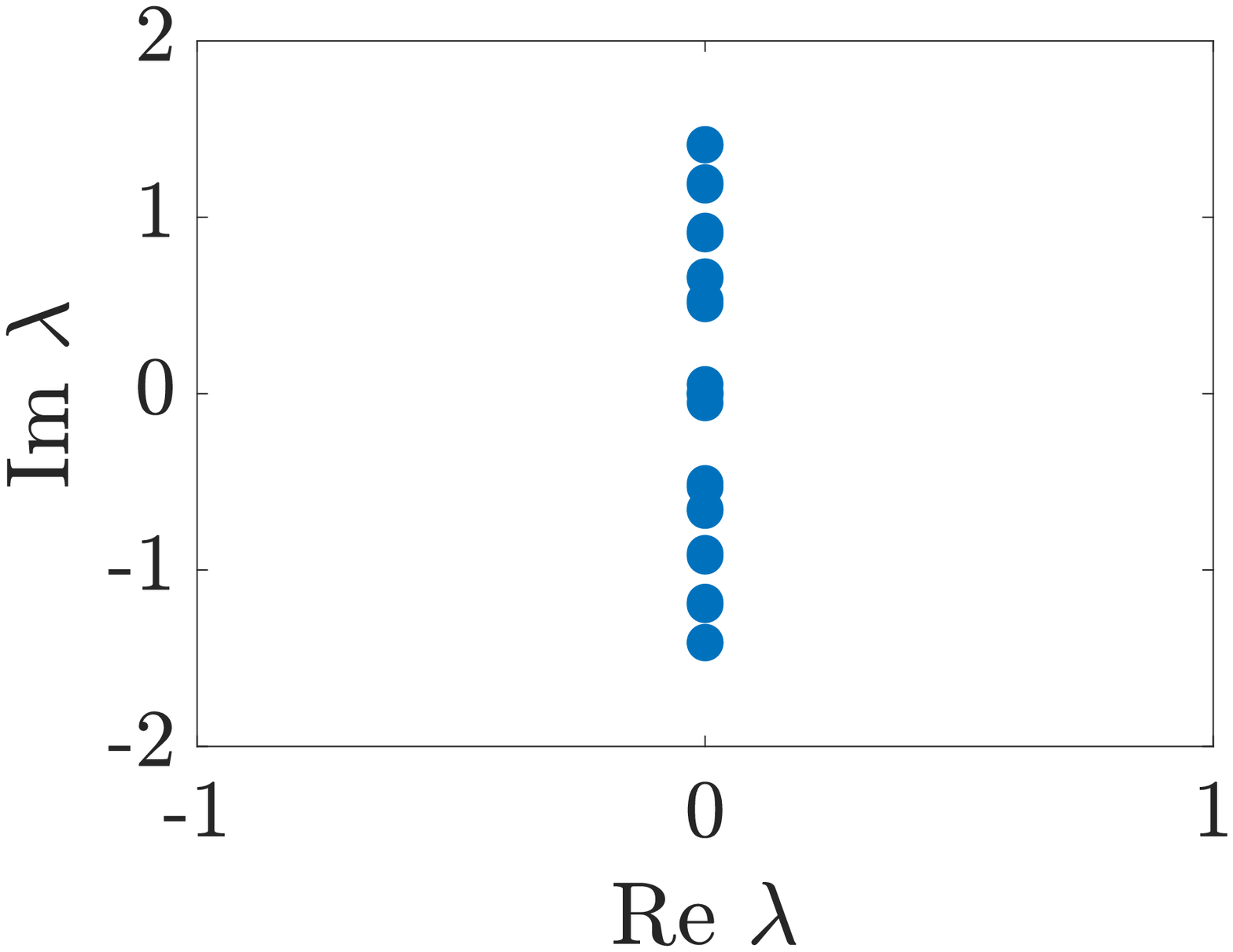}
\includegraphics[width=4cm]{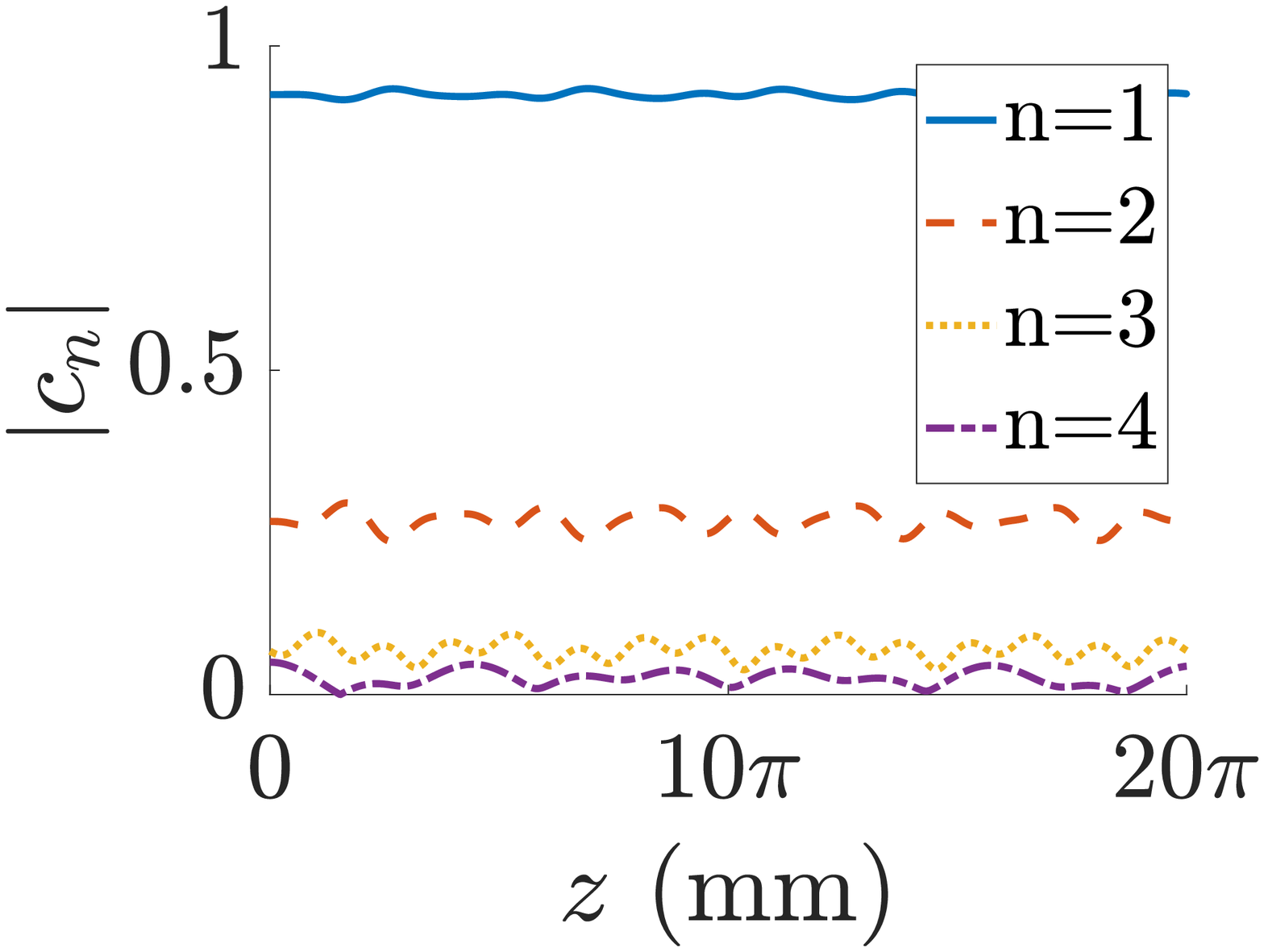}
\end{tabular}
\end{center}
\caption{Spectrum of linearization of \cref{eq:twist1} about symmetric double pulse solution with $N=12$ and $\phi = \pi/6$ (left panel). Amplitude $|c_n|$ for first four nodes versus $z$ for evolution of perturbation of this solution using fourth order Runge-Kutta scheme (right panel). $\omega = 1\text{ mm}^{-1}$, $k=0.25\text{ mm}^{-1}$.}
\label{fig:dpstab}
\end{figure}

Numerical spectral computations, as well as numerical evolution of perturbations of these solutions, suggest that these double pulse solutions are neutrally stable (\cref{fig:dpstab}).

\section{Conclusions}

In this paper, we have demonstrated the existence of standing wave solutions to a system of equations modeling light propagation in a twisted multi-core fiber in the setting of no gain or loss at the individual sites. Our theoretical results extend previous work and add understanding on stability properties. It is both intriguing and fascinating that the mathematical tool used here(continuation) to build exact solutions discovers, in a natural way, a physical phenomenon (AB suppression). The mathematical approach reveals the role of symmetries, phase relations and nonlinearity; the last one is evident in what is used as the starting ($k=0$) solution for the continuation method. We find specifically that if the twist parameter $\phi$ and the number of waveguides $N$ are related by $\phi = \pi/N$, then standing wave solutions exist which are a manifestation of the optical Aharonov-Bohm suppression, i.e. there is a node which is completely dark for all time. These solutions exist for both $N$ even and $N$ odd, and are all neutrally stable. While we emphasize the theory here, suitable parameters and powers for experimental realizations suggested in, for example, \cite[Figure 3]{Parto2019} should apply for a range of values shown here (e.g. $N = 6-10$). For future research, it would be interesting to investigate whether such standing waves exist for twisted optical fibers in more complicated geometries such as multiple concentric rings or Lieb lattices. We could also systematically study multi-pulse solutions, as well as investigate the existence and stability of breathers, which are localized, periodic structures that are not standing waves. (See \cite{Lumer2013} for examples of breather solutions in honeycomb lattices). We could also apply the techniques used here to the $\mathcal{PT}$-symmetric system with symmetric gain and loss, which is studied in \cite{castro2016}. Finally, since these standing wave solutions are neutrally stable, it would be interesting to see if they could be created experimentally in twisted multi-core fibers.

\begin{acknowledgments}
This material is based upon work supported by the U.S. National Science Foundation under the RTG grant DMS-1840260 (R.P. and A.A.) and DMS-1909559 (AA). The authors would also like to thank P.G. Kevrekidis for his helpful comments and suggestions for numerical simulations.
\end{acknowledgments}

\bibliographystyle{apsrev4-2}
\bibliography{twist.bib}

\end{document}